\documentclass[12pt,a4paper]{amsart}

\usepackage{tikz}
\usetikzlibrary{arrows}
\usetikzlibrary{calc}

\usepackage{graphicx}
\usepackage{amsmath}
\usepackage{verbatim}
\usepackage{colonequals}

\theoremstyle{plain}
\newtheorem{thm}{Theorem}[section]

\newtoks\prt
\newtheorem{proclaim}[thm]{\the\prt}
\theoremstyle{definition}

\def\eqn#1$$#2$${\begin{equation}\label#1#2\end{equation}}

\numberwithin{equation}{section}

\def\adj{\operatorname{adj}}

\def\cof{\operatorname{cof}}
\def\dadj{\mathcal{ADJ}\, }

\def\djac{\mathcal{J}}
\def\J{\mathcal{J}}

\def\diam{\operatorname{diam}}
\def\dist{\operatorname{dist}}

\def\epsilon{\varepsilon}

\def\er{\mathbb R}
\def\es{\mathbb S}

\def\haus{\mathcal{H}}

\def\loc{\operatorname{loc}}
\def\M{\mathcal{M}}
\def\mir2{\mathcal{L}^2}

\def\phi{\varphi}

\def\r2{\er^2}
\def\sgn{\operatorname{sgn}}

\def\spt{\operatorname{supp}}
\def\supp{\operatorname{supp}}
\def\testers{C_0^{\infty}}

\def\rn{\mathbb R^n}
\def\spt{\operatorname{supp}}
\def\sgn{\operatorname{sgn}}

\def\zet{\mathbb Z}

\newcommand{\norm}[1]{\left\Arrowvert #1 \right\Arrowvert}
\newcommand{\abs}[1]{\left\vert #1 \right\vert}
\newcommand{\inv}{^{-1}}
\newcommand{\eps}{\varepsilon}

\newtoks\by
\newtoks\paper
\newtoks\book
\newtoks\jour
\newtoks\yr
\newtoks\pages
\newtoks\vol
\newtoks\publ
\def\ota{{\hbox\vol{???}}}
\def\cLear{\by=\ota\paper=\ota\book=\ota\jour=\ota\yr=\ota
  \pages=\ota\vol=\ota\publ=\ota}
\def\endpaper{\the\by, {\the\paper},
  \textit{\the\jour} \textbf{\the\vol} (\the\yr), \the\pages.\cLear}
\def\endbook{\the\by, \textit{\the\book}, \the\publ.\cLear}
\def\endprep{\the\by, \textit{\the\paper}, \the\jour.\cLear}
\def\endyearprep{\the\by, \textit{\the\paper}, \the\jour, (\the\yr).\cLear}
\def\name#1#2{#2 #1}

\title
{Weak regularity of the inverse under minimal assumptions}
\author[S. Hencl]{Stanislav Hencl}
\address{Department of Mathematical Analysis, Charles University, Sokolovsk\'{a}
  83, 186 00
  Prague 8, Czech Republic}
\email{hencl@karlin.mff.cuni.cz}

\author[A. Kauranen]{Aapo Kauranen}
\address{Department of Mathematical Analysis, Sokolovska 83, Praha 8, 186 75, Charles University in Prague \and
  Departament de Matem\`atiques, Universitat Aut\`onoma de Barcelona, 08193, Bellaterra (Barcelona), Spain}
\email{aapo.p.kauranen@gmail.com}

\author[R. Luisto]{Rami Luisto}
\address{Department of Mathematical Analysis, Charles University, Sokolovsk\'{a}
  83, 186 00
  Prague 8, Czech Republic}
\email{rami.luisto@gmail.com}

\thanks{SH and AK  were supported in part by the ERC CZ grant LL1203 of the
Czech Ministry  of Education. SH  was supported in part by the grant
GA\v{C}R P201/18-07996S.
AK acknowledges financial support from the Spanish Ministry of
Economy and Competitiveness, through the ``Mar\'{\i}a de Maeztu'' Programme for Units of
Excellence in R\&D (MDM-2014- 0445).
The last author was supported by a grant of the
Finnish Academy of Science and Letters.}

\begin{document}

\begin{abstract}
  Let $\Omega\subset\er^3$ be a domain and let $f\in BV_{\loc}(\Omega,\er^3)$ be a
  homeomorphism such that its distributional adjugate is a finite Radon measure.
  We show that its inverse has bounded variation $f^{-1}\in BV_{\loc}$.
 The condition that the distributional adjugate is finite measure is not only sufficient but also necessary for the weak regularity of the inverse.
\end{abstract}

\maketitle

\section{Introduction}
Suppose that $\Omega\subset \rn$ is an open set
and let $f \colon \Omega\to f(\Omega)\subset \rn$ be a homeomorphism.
In this paper we address the issue of
the weak regularity of $f^{-1}$ under regularity assumptions on $f$.

The classical inverse function theorem states that the inverse of a
$C^1$-smooth homeomorphism $f$ is again a $C^1$-smooth homeomorphism,
under the assumption that the Jacobian $J_f$ is strictly positive.
In this paper we address the question whether the inverse of a
Sobolev or $BV$-homeomorphism
is a $BV$ function or even a Sobolev function. This problem is of
particular importance as Sobolev and $BV$ spaces are commonly used as initial
spaces for existence problems in PDE's and the calculus of variations.
For instance, elasticity is a typical field where
both invertibility problems and Sobolev (or $BV$) regularity issues
are relevant (see e.g.\ \cite{Ball}, \cite{CiaNec} and \cite{MST}).

The problem of the weak regularity of the inverse has attracted a big attention
in the past decade. It started with the result of \cite{HK} and \cite{HKO} where
it was shown that for homeomorphisms  in dimension $n=2$ we have
$$
\begin{aligned}
  &\bigl(f\in W^{1,1}_{\loc}(\Omega,\er^2)\text{ is a mapping of finite distortion
  }\Rightarrow
  f^{-1}\in W^{1,1}_{\loc}(f(\Omega),\er^2)\bigl)\\
  &\text{ and }\bigl(f\in BV_{\loc}(\Omega,\er^2)\Rightarrow f^{-1}\in
  BV_{\loc}(f(\Omega),\er^2)\bigr).\\
\end{aligned}
$$
This result has been generalized to $\rn$ in \cite{CHM} where it was shown
$$
\begin{aligned}
  &\bigl(f\in W^{1,n-1}_{\loc}(\Omega,\rn)\text{ is a mapping of finite distortion
  }\Rightarrow
  f^{-1}\in W^{1,1}_{\loc}(f(\Omega),\rn)\bigl)\\
  &\text{ and }\bigl(f\in W^{1,n-1}_{\loc}(\Omega,\rn)\Rightarrow f^{-1}\in
  BV_{\loc}(f(\Omega),\rn)\bigr).\\
\end{aligned}
$$

It is natural to study the sharpness of the assumption $f\in W^{1,n-1}$. It was
shown in \cite{H2} that the assumption $f\in W^{1,n-1-\epsilon}$ (or any weaker
Orlicz-Sobolev assumption) is not sufficient in general. Furthermore, by results
of \cite{DS} we know that for $f\in W^{1,n-1}$ we have not only $f^{-1}\in BV$
but also the total variation of the inverse satisfies
\eqn{identity}
$$
|Df^{-1}|(f(\Omega))=\int_{\Omega}|\adj Df(x)|\; dx
$$
where $\adj Df$ denotes the adjugate matrix, i.e. the matrix of $(n-1)\times(n-1)$
subdeterminants. However, for $n\geq 3$ it is possible to construct a $W^{1,1}$
homeomorphism with $\adj Df\in L^1$ such that $f^{-1}\notin BV$ (see \cite{H2})
so the pointwise adjugate does not carry enough information about the regularity
of the inverse. The main trouble in the example from \cite{H2} is that the pointwise
adjugate does not capture some singular behavior on the set of measure zero.

On the other hand, the results of \cite{CHM} are not perfect as they cannot be applied to even very
simple mappings. Let $c(x)$ denote the usual Cantor ternary function, then
$h(x)=x+c(x)$ is $BV$ homeomorphism and its inverse $g=h^{-1}$ is even Lipschitz.
It is easy to check that
$$
f(x,y,z)=[h(x),y,z]
$$
is a $BV$ homeomorphism and its inverse $f^{-1}(x,y,z)=[g(x),y,z]$ is Lipschitz,
but the results of \cite{CHM} cannot be applied as $f$ is not Sobolev.  In this
paper we obtain a new result in dimension $n=3$ about the regularity of the inverse which
generalizes the result of \cite{CHM} and can be applied to the above mapping.

It is well-known that in models of Nonlinear Elasticity and in Geometric Function
Theory the usual pointwise Jacobian does not carry enough information about the
mapping and it is necessary to work with the distributional Jacobian, see e.g.\
\cite{Ball}, \cite{CL}, \cite{HKbook}, \cite{IKO}, \cite{IM} and \cite{M}.
This distributional Jacobian captures the behavior on zero measure sets and
can be used to model for example cavitations of the mapping, see e.g.\
\cite{HeMo1}, \cite{HeMo2}, \cite{MS} and \cite{MST}.
In the same spirit we introduce
the notion of the distributional adjugate $\dadj Df$ (see Definition \ref{def} below)
and we show that the right assumption for the regularity of the inverse is that
$\dadj Df\in \M(\Omega,\er^{3\times 3})$, where $\M(\Omega)$ denotes finite Radon
measures on $\Omega$.

Further we need to add the technical assumption that Lebesgue area (see \eqref{lebarea} below) of image of a.e. hyperplane parallel to coordinate axes is finite. Let us recall that the Hausdorff measure is always bigger (see section \ref{sec:Areas}) 
$$
A(f,\Omega\cap\{x\in\er^3:x_j=t\})\leq \haus^2\bigl(f(\Omega\cap\{x\in\er^3:x_j=t\})\bigr),\ j=1,2,3
$$ 
so it is enough to assume the finiteness of Hausdorff measure of the image. 
 Our main result is the following.
\prt{Theorem}
\begin{proclaim}\label{regularity}
  Let $\Omega\subset\er^3$ be a domain and $f\in BV_{\loc}(\Omega,\er^3)$ be a
  homeomorphism such that $\dadj Df\in \M(\Omega,\er^{3\times 3})$  and assume further that for a.e. $t$ we have
	\eqn{stupid}
	$$
	A(f,\Omega\cap\{x\in\er^3:x_j=t\})<\infty\text{ for }j=1,2,3. 
	$$
	Then $f^{-1}\in
  BV_{\loc}(f(\Omega),\er^3)$.

  If we moreover know that the image of the measure $f(\dadj Df)$ is absolutely
  continuous with respect to Lebesgue measure, then $f^{-1}\in
  W_{\loc}^{1,1}(f(\Omega),\er^3)$.
\end{proclaim}
It would be very interesting to see if the assumption \eqref{stupid} can be removed. A similar extra assumption was assumed in \cite[Theorem 14]{DeLellis}.

To show that our result generalizes the aforementioned result of \cite{CHM}
we notice first that for homeomorphisms in $W^{1,n-1}$ the
distributional adjugate $\dadj Df$ is equal to the pointwise adjugate $\adj Df$
(see \cite[Proposition 2.10]{HKbook}).  The main part of the proof in \cite{CHM} was to show that
$f$ maps $\haus^{n-1}$ null sets on almost every hyperplane to
$\haus^{n-1}$ null sets. This implies  \eqref{stupid} for $W^{1,n-1}$ -homeomorphisms.
This property of null sets on hyperplanes may fail in our setting so our proof
is more subtle and we have to use delicate tools of Geometric Measure Theory.
Moreover, our assumptions are
not only sufficient but also necessary for the weak regularity of the inverse.

\prt{Theorem}
\begin{proclaim}\label{reverse}
  Let $\Omega\subset\er^3$ be a domain and $f\in BV(\Omega,\er^3)$ be a
  homeomorphism such that $f^{-1}\in BV(f(\Omega),\er^3)$.
  Then $\dadj Df\in \M(\Omega,\er^{3\times 3})$ and for a.e. $t$ we have
	$$
	A(f,\Omega\cap\{x\in\er^3:x_j=t\})<\infty\text{ for }j=1,2,3. 
	$$
\end{proclaim}

In an upcoming article \cite{HKMdraft} these results are further refined. There the total variation of distributional adjugate is shown to \emph{equal} the total variation of the derivative of the inverse mapping. Moreover,  some simple ways of verifying the key assumption $\dadj Df\in \M(\Omega,\er^{3\times 3})$  are presented there.

Now we give the formal definition of the distributional adjugate. Without loss of
generality we can assume that $\Omega=(0,1)^3$ as all statements are local.

\begin{figure}[h!]
  \centering
  \resizebox{\textwidth}{!}{
    \begin{tikzpicture}
      % The left and right blocks are both mirrored w.r.t. to y-axis,
      % but the function arrows are not. This is due to lack of planning,
      % but makes altering the picture a nice little puzzle.

      % Changing this scale changes the font size of mathematical symbols.
      % Smaller number means bigger symbols.
      \begin{scope}[scale=1.2]

        % Left block
        \begin{scope}[shift={(-2,0)},xscale=-1]

          % Constants controlling perspective.
          \def\vert{0.707}
          \def\hori{0.707}

          % Height of the restriction slab.
          \def\t{1.1}

          % Front face of the cube
          \coordinate (A) at (0,2);
          \coordinate (B) at (2,2);
          \coordinate (C) at (0,0);
          \coordinate (D) at (2,0);

          % Back face of the cube
          \coordinate (E) at ($(A)+(\vert,\hori)$);
          \coordinate (F) at ($(B)+(\vert,\hori)$);
          \coordinate (G) at ($(C)+(\vert,\hori)$);
          \coordinate (H) at ($(D)+(\vert,\hori)$);

          % Coordinates of the gray slab.
          \foreach \P in {C,D,G,H}{
          \coordinate (t\P) at ($(\P) + (0,\t)$);
          }

          % Drawing dashed lines in the back for cube
          \foreach \P in {H,C,E} {
            \draw[dashed] (G) -- (\P);
          }

          % Drawing the gray slab
          \draw[fill, opacity=0.3, black!50] (tC)--(tD)--(tH)--(tG)--cycle;
          \draw[thick] (tC)--(tD)--(tH);
          \draw[thick, dashed] (tH)--(tG)--(tC);

          % Drawing the cube
          \draw[] (A)--(B)--(D)--(C)--cycle;
          \draw[] (A)--(B)--(F)--(E)--cycle;
          \draw[] (B)--(F)--(H)--(D)--cycle;

          % Marking the slab
          \node[left] at (tH) {$t$};

          % End of left block
        \end{scope}

      % The function arrows for $f$.
      \draw[->, thick] (-1.5,1.1) to [out=10,in=170] (2,1.5);
      \draw[->, thick] (-1.5,1.6) to [out=60,in=180] (0.8,3);
      \draw[->, thick] (-1.5,0.6) to [out=-60,in=-180] (2,-1);

      \node[] at (0,1.1) {$f^t_3$};
      \node[above] at (-0.2,2.1) {$f^t_{3,2}$};
      \node[above] at (0,-0.6) {$f^t_{3,3}$};
      \node[right] at (7,1) {$f^t_{3,1}$};

      % Right Block
      \begin{scope}[shift={(5,0)},xscale=-1]

        % Same perspective constant as in left block
        \def\vert{0.707}
        \def\hori{0.707}

        % Height of gray slab
        \def\t{0.9}

        % This allows to change the place of the image of the cube.
        % Making it non-zero will break part of the picture since
        % some parts have literal constants. Sorry.
        \coordinate (S) at (0,0);

        % Front face of image cube
        \coordinate (A) at ($(0,2) + (S)$);
        \coordinate (B) at ($(2,2) + (S)$);
        \coordinate (C) at ($(0,0) + (S)$);
        \coordinate (D) at ($(2,0) + (S)$);

        % Back face of image cube
        \coordinate (E) at ($(A)+(\vert,\hori)$);
        \coordinate (F) at ($(B)+(\vert,\hori)$);
        \coordinate (G) at ($(C)+(\vert,\hori)$);
        \coordinate (H) at ($(D)+(\vert,\hori)$);

        % Image of gray slab.
        \foreach \P in {C,D,G,H}{
          \coordinate (t\P) at ($(\P) + (0,\t)$);
        }

        % Drawing the image of the cube.
        % Hand-picked constants, fiddle if you like.

        % Bottom
        \draw[ dashed]
        (C) to [out = 40, in=220]
        (G) to [out = -15, in=150]
        (H);
        \draw[]
        (H) to [out = 200, in=55]
        (D) to [out = 190, in=20]
        (C);

        % Top
        \draw[]
        (A) to [out = 30, in=180]
        (E) to [out = 20, in=200]
        (F) to [out = 210, in=45]
        (B) to [out = 130, in=-10]
        (A);

        % Vertical - front left
        \draw
        (C)  to [out=95, in=260]
        (tC) to [out=80, in=260]
        (A);

        % Vertical - front right
        \draw
        (D)  to [out=95, in=260]
        (tD) to [out=80, in=260]
        (B);

        % Vertical - back right
        \draw
        (H)  to [out=95, in=260]
        (tH) to [out=80, in=260]
        (F);

        % Vertical - back left
        \draw[dashed]
        (G)  to [out=95, in=260]
        (tG) to [out=80, in=260]
        (E);

        % Drawing of cube finished

        % x-direction projection
        \coordinate (LL) at (-2,0);
        \coordinate (LU) at (-2,2);
        \coordinate (LLo) at ($(LL) + (\vert,\hori)$);
        \coordinate (LUo) at ($(LU) + (\vert,\hori)$);

        % x-direction Boundary
        \draw[thin,dashed]
        (LL)  to [out=95,in=280]
        (LU)  to [out =35, in = 210]
        (LUo) to [out =260, in = 80]
        (LLo) to [out =215, in = 50]
        (LL);

        % x-direction projection lines
        \draw[thin,dotted] (C) -- (LL);
        \draw[thin,dotted] (A) -- (LU);
        \draw[thin,dotted] (E) -- (LUo);
        \draw[thin,dotted] (G) -- (LLo);

        % x-direction Image of gray slab
        \draw[fill,black!30,opacity=0.3]
        ($(tC) - (1.95,0)$) to [out=80,in=280]
        ($(tC) - (1.95,-0.1)$) to [out=30, in= 200]
        ($(tG) - (2,0)$) to [out=260,in=100]
        ($(tG) - (2,0.2)$) to [out=200, in=15]
        ($(tC) - (1.95,0)$);
        \draw
        ($(tC) - (1.95,0)$) to [out=80,in=280]
        ($(tC) - (1.95,-0.1)$) to [out=30, in= 200]
        ($(tG) - (2,0)$) to [out=260,in=100]
        ($(tG) - (2,0.2)$) to [out=200, in=15]
        ($(tC) - (1.95,0)$);

        % z-directional projection
        \coordinate (BL) at (0,-1.5);
        \coordinate (BR) at (2,-1.5);
        \coordinate (BLo) at ($(BL) + (\vert,\hori)$);
        \coordinate (BRo) at ($(BR) + (\vert,\hori)$);

        % z-direction boundary box
        \draw[thin,dashed]
        (BL)  to [out=5,in=190]
        (BR)  to [out=65,in=250]
        (BRo) to [out=170,in=-10]
        (BLo) to [out=40,in=190]
        (BL);

        % z-direction image of gray slab
        \draw[fill,black!30,opacity=0.3]
        ($(BL)+(0.1,-0.1)$)  to [out=-5,in=170]
        ($(BR)+(-0.1,-0.1)$)  to [out=75,in=270]
        ($(BRo)+(0.05,-0.1)$) to [out=160,in=-5]
        ($(BLo)+(0.1,0)$) to [out=250,in=100]
        ($(BL)+(0.1,-0.1)$);
        \draw
        ($(BL)+(0.1,-0.1)$)  to [out=-5,in=170]
        ($(BR)+(-0.1,-0.1)$)  to [out=75,in=270]
        ($(BRo)+(0.05,-0.1)$) to [out=160,in=-5]
        ($(BLo)+(0.1,0)$) to [out=250,in=100]
        ($(BL)+(0.1,-0.1)$);

        % z-direction projection lines
        \draw[thin,dotted] (C) -- (BL);
        \draw[thin,dotted] (D) -- (BR);
        \draw[thin,dotted] (H) -- (BRo);
        \draw[thin,dotted] (G) -- (BLo);

        % y-directional projection
        \foreach \P in {E,F,G,H} {
          \coordinate (T\P) at ($(\P) + (2*\vert,2*\hori)$);
        }
        \foreach \P in {E,F,G,H} {
          \draw[thin,dotted] (\P)--(T\P);
        }

        % y-direction bounding box
        \draw[thin,dashed]
        (TE) to [out=5, in=170]
        (TF) to [out = 250, in=80]
        (TH) to [out = 190, in=-10]
        (TG) to [out = 95, in=260]
        (TE);

        % y-directional image of gray slab
        \draw[thick]
        ($(tG) + (2*\hori,2*\vert)$)
        to [ out= -10, in = 150]
        ($(tH) + (2*\hori,2*\vert)$) ;

        % Main image, drawing the image of the gray slab
        \draw[fill, black!30, opacity=0.3]
        (tC) to [out = 70, in=200]
        (tG) to [out = -20, in=120]
        (tH) to [out = 210, in=45]
        (tD) to [out = 100, in=-10]
        (tC);
        \draw[thick, dashed]
        (tC) to [out = 70, in=200]
        (tG) to [out = -20, in=120]
        (tH);
        \draw[thick]
        (tH) to [out = 210, in=45]
        (tD) to [out = 100, in=-10]
        (tC);

        % Right block ends
      \end{scope}

      % Picture ends
      \end{scope}
    \end{tikzpicture}
  }
  \label{fig:NineMappings}
  \caption{The definition of the nine 2-dimensional restrictions $f^t_{k,j}$ of a mapping $f \colon (0,1)^3 \to \mathbb{R}^3$
  in Definition \ref{slicedef}.}
\end{figure}
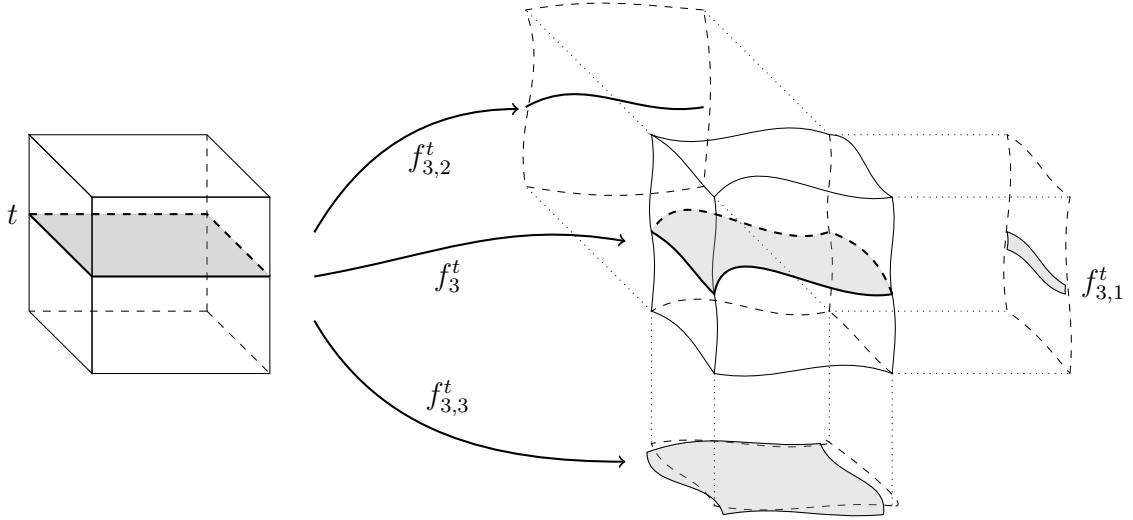

\prt{Definition}
\begin{proclaim}\label{slicedef}
  Let $f=(f_1,f_2,f_3) \colon (0,1)^3\to\er^3$ be a homeomorphism in $BV$. For $t\in (0,1)$ we define
  $$
  f^{t}_1(x)=f(t,x_2,x_3),\ f^{t}_2(x)=f(x_1,t,x_3)\text{ and }f^{t}_3(x)=f(x_1,x_2,t).
  $$
  We can split these mappings into $9$ mappings from $(0,1)^{2}\to\er^{2}$ using
  its coordinate functions. Given $k,j\in\{1,2,3\}$ choose $a,b\in\{1,2,3\}\setminus\{j\}$ with $a<b$ and define
  $$
  f^{t}_{k,j}(x)=\bigl[(f^{t}_{k})_a(x),(f^{t}_{k})_{b}(x)\bigr],
  $$
  see Figure \ref{fig:NineMappings}. For example,
  $$
  f^{t}_{1,1}(x_2,x_3)=\bigl[(f^{t}_{1})_2(x),(f^{t}_{1})_{3}(x)\bigr]=\bigl[f_2(t,x_2,x_3),f_{3}(t,x_2,x_3)\bigr].
  $$
  \end{proclaim}

  Now we recall the definition of distributional Jacobian and, using it, define
the distributional adjugate.

  \prt{Definition}
  \begin{proclaim}\label{def}
  Let $f$ be as in  Definition \ref{slicedef}.
  For  mappings $f^t_{k,j}$ we consider the usual distributional Jacobian (see
e.g.\
  \cite[Section 2.2]{HKbook}), i.e. the distribution
  $$
  \djac_{f^{t}_{k,j}}( \varphi) = -\int_\Omega (f^{t}_{k,j})_1(x)
  J\bigl(\varphi, (f^{t}_{k,j})_2\bigr)(x)\; dx \qquad \text{ for
    all } \varphi \in \testers(\Omega).
  $$
  This distribution is well-defined for homeomorphism in $W^{1,1}$. It
  is also well-defined for homeomorphism in $BV$ for $n=3$, we just consider the
  integral with respect to corresponding measure $d(\partial_l f^{t}_{k,j})_{2}(x)$
  instead of $(\partial_l f^{t}_{k,j})_2(x)\; dx$ and for example we define
  $$
  \begin{aligned}
    \djac_{f^{t}_{1,1}}( \varphi)=&-\int_\Omega f_2(t,x_2,x_3)
    \frac{\partial\varphi(x_2,x_3)}{\partial x_2}d\Bigl(\frac{\partial f_3(t,x_2,x_3)}{\partial x_3}\Bigr)\\
    &+\int_\Omega f_2(t,x_2,x_3)
    \frac{\partial\varphi(x_2,x_3)}{\partial x_3}d\Bigl(\frac{\partial f_3(t,x_2,x_3)}{\partial x_2}\Bigr).\\
  \end{aligned}
  $$

  Assume that these $3\times 3$ distributions $\djac_{f^{t}_{i,j}}$ are measures
  for a.e.\ $t\in(0,1)$ and for measurable $A\subset(0,1)^3$ we set
  $$
  (\dadj Df)_{k,j}(A)=\int_0^1 \djac_{f^{t}_{k,j}}\left(A\cap\{x_k=t\}\right)\; dt.
  $$
  We say that $\dadj Df\in \M(\Omega,\er^{3\times 3})$ if the distributions
  $\djac_{f^{t}_{k,j}}$ are measures for a.e.\ $t\in(0,1)$ and $(\dadj Df)_{k,j}\in
  \M(\Omega)$ for every $i,j\in\{1,2,3\}$.
\end{proclaim}

A priori it seems that the definition of the distributional adjugate
is dependent on the choice of coordinates. This turns out not to be the
case and we discuss this further in Section \ref{sec:ADJ-CoordinateDependence}.

\section{Preliminaries}

Total variation of the measure $\mu$ is the measure $|\mu|$ such that
$$
|\mu|(A)\colonequals  \sup\left\{\int_{\rn} \phi \; d\mu:\ \phi\in C_0(A),\
\|\phi\|_{\infty}\leq 1\right\}\text{ for all open sets }A\subset\rn .
$$

For a domain $\Omega \subset \rn$ we denote by $C^\infty_0(\Omega)$ those smooth
functions $\varphi$ whose support is compactly contained in $\Omega$,
i.e.\ $\supp \varphi \subset \subset \Omega$.

\subsection{Mollification}

We will need to approximate continuous $BV$ mappings with smooth
maps. To this end we recall here the
basic definitions of convolution and mollifiers
for the reader's convenience; for a more detailed treatise on
the basics and connections to $BV$ mappings we refer to \cite[Sections 2.1 and 3.1]{AFP}.

A family of mappings $(\rho_\eps) \in C^\infty(\er^n, \er)$
is called \emph{a family of mollifiers} if for all $\eps>0$ we have
$\rho_\eps(x) = \eps^{-n} \rho(x / \eps)$, where $\rho \in C^\infty(\er^n, \er)$
is a non-negative mapping satisfying
$\supp \rho \subset B({0},1)$,
$\rho(-x) = \rho(x)$ and
$\int_{\rn} \rho = 1$.
We will sometimes use a \emph{sequence} of mollifiers $(\rho_j)$,
in which case we tacitly assume that there is a family of mollifiers $(\tilde\rho_\eps)$
from which we extract the sequence $(\rho_j)$ by setting $\rho_j \colonequals \tilde\rho_{\frac 1 j }$.

For $\Omega \subset \rn$ and any two functions
$f \colon \Omega \to \er^m$, $g \colon \Omega \to \er$
we set their \emph{convolution} to be
\begin{align*}
  (f \ast g) \colon \rn \to \er^m,
  \quad
  (f \ast g)(x)
  = \int_{\Omega} f(y) g(x-y) \; dy,
\end{align*}
whenever the integral exists.
Likewise for a $m$-valued Radon measure $\mu$ defined on $\Omega$
and a function $g \colon \Omega \to \er$
we define their convolution as
\begin{align*}
  (\mu \ast f) \colon \rn \to \er^m,
  \quad
  (\mu \ast f)(x)
  = \int_{\Omega} f(x-y) \; d\mu(y)
\end{align*}
whenever the integral exists.

For a function $f \colon \Omega \to \er^m$
or a Radon measure $\mu$ defined on $\Omega$
we define their \emph{(family of) mollifications}
to be the families
$(f_\eps) \colonequals (f \ast \rho_\eps)$
and
$(\mu_\eps)\colonequals (\mu \ast \rho_\eps)$,
respectively, where $(\rho_\eps)$
is a family of mollifiers. Similarly we define the
sequence of mollifications as
$(f_j) \colonequals (f \ast \rho_j)$
and
$(\mu_j) \colonequals (\mu \ast \rho_j)$.
For our purposes the exact family of mollifiers
does not matter, so we tacitly assume that some such family
has been given whenever we use mollifications.

\subsection{Topological degree}
\label{sec:TopologicalDegree}

For $\Omega \subset \rn$ and a given smooth map $f \colon \Omega \to \rn$
we define the \textit{topological degree} as
\begin{align*}
  \deg(f,\Omega, y_{0}) = \sum_{x \in \Omega \cap f\inv\{y_0\}} \sgn(J_{f}(x))
\end{align*}
if $J_{f}(x) \neq 0$ for each $x \in f^{-1}\{y_{0}\}$. This definition can be extended to arbitrary continuous mappings and each point $y_0\notin f(\partial\Omega)$, see e.g.\ \cite[Section 1.2]{FG} or \cite[Chapter 3.2]{HKbook}.
For our purposes the following property of the topological degree is crucial;
see \cite[Definition 1.18]{FG}.
\prt{Lemma}
\begin{proclaim}
  \label{lemma:TopologicalDegree}
  Let $\Omega \subset \rn$ be a domain,
  $f \colon \Omega \to \rn$ a continuous function and
  $U$ a domain with $\overline{U} \subset \Omega$.
  Then for any point $y_0 \in \rn \setminus f (\partial U)$
  and any continuous mapping $g \colon \Omega \to \rn$ with
  \begin{align*}
    \| f - g \|_\infty
    \leq \operatorname{dist}\left( y_0 , f(\partial U) \right),
  \end{align*}
  we have $\deg(f, \Omega, y_0) = \deg(g, \Omega, y_0)$.
\end{proclaim}

We will also need some classical results concerning
the dependence of the degree on the domain.
The following result is from \cite[Theorem 2.7]{FG}.
\prt{Lemma}
\begin{proclaim}
  \label{lemma:TopologicalDegree2}
  Let $\Omega \subset \rn$ be a domain,
  $f \colon \Omega \to \rn$ a continuous function and
  $U$ a domain with $\overline{U} \subset \Omega$.
  \begin{enumerate}
  \item \emph{(Domain decomposition property)} For any domain $D \subset U$ with a decomposition
  $D = \cup_i D_i$ into open disjoint sets, and
  a point $ p \notin f (\partial D)$, we have
  \begin{align*}
    \deg(f,D,p)
    = \sum_i \deg(f,D_i,p).
  \end{align*}

  \item \emph{(Excision property)} For a compact set $K \subset \overline U$ and
    a point $p \notin f(K \cup \partial U)$  we have
    $\deg(f,U,p) = \deg(f,U\setminus K,p)$.
  \end{enumerate}
\end{proclaim}

The topological degree agrees with the Brouwer
degree for continuous mappings, which in turn equals
the \emph{winding number} in the plane.
The winding number is an integer expressing how many times the path $\beta_f \colonequals f (\partial D)$
circles the point $p$; indeed, the winding number equals the topological index
of the mapping $\frac{\beta_f - p}{|\beta_f - p|} \colon \es^1 \to \es^1$.
We refer to \cite[Section 2.5]{FG} for discussion of the winding number in the
setting of holomorphic planar mappings.

\subsection{Hausdorff measure}
For $A\subset\rn$ we use the classical definition of the Hausdorff measure (see e.g.\ \cite{Fe})
$$
\haus^k(A)=\lim_{\delta\to 0+}\haus^k_{\delta}(A),
$$
where
$$
\haus^k_{\delta}(A)=\inf\Bigl\{\sum_i \diam^k A_i:\ A\subset\bigcup_i A_i,\ \diam A_i\leq \delta\Bigr\}.
$$
The important ingredient of our proof is the Gustin boxing inequality \cite{G} which states that for each compact set $K\subset \rn$ we have
\eqn{gustin}
$$
\haus^{n-1}_{\infty}(K)\leq C_n \haus^{n-1}(\partial K).
$$

\subsection{On various areas}
\label{sec:Areas}

Besides homeomorphisms in three dimensions we work with a continuous mappings $g\colon[0,1]^2\to\er^3.$ A central object is the Hausdorff measure of the image $\haus^2(g((0,1)^2)),$ but we need to use other finer notions of area. The results of this subsection can be found in the book of Cesari \cite{C} and they also follow by some results of Federer, see e.g.\ \cite[(13) on page 93]{Fe2} and references given there.

First we define the Lebesgue area (see \cite[3.1]{C}). Let $L$ be an affine
mapping, then for any triangle $\Delta$ the area of $L(\Delta)$ is defined in
the natural way.   For a piecewise linear mapping $h:[0,1]^2\rightarrow \er^3$
we define the Lebesgue area $A(h, [0,1]^2)$ to be the sum of the areas of the triangles
of some triangulation where $h$ is linear in each of these triangles. We define
\begin{align}\label{lebarea}
  A(g, [0,1]^2)
  \colonequals \inf
  \left\{
  \lim_{k \to \infty}A(g_k, [0,1]^2) : (g_k) \in \mathcal{PH}(g)
  \right\},
\end{align}
where $\mathcal{PH}(g)$ is the collection of all sequences of polyhedral surfaces converging uniformly to $g$.

Next we define coordinate mappings $g_{j}\colon[0,1]^2\to\er^2$ as
$$
g_{1}(x)=[g_2(x),g_3(x)],\ g_{2}(x)=[g_1(x),g_3(x)]\text{ and }g_{3}(x)=[g_1(x),g_2(x)].
$$

Finally we define (see \cite[9.1]{C})
\eqn{recalldef}
$$
V(g_{j},[0,1]^2)\colonequals  \sup_S\left\{\sum_{\pi\in S}\int_{\er^2}
\bigl|\deg(g_{j},\pi,y)\bigr|\; dy\right\},
$$
where $S$ is any finite system of nonoverlapping simple open polygonal regions in
$[0,1]^2$ and
$\deg(g_{j},y,A)$ denotes the topological degree of mapping.

We need the following characterization of the Lebesgue area (see \cite[18.10
and 12.8.$(ii)$]{C}) which holds for any continuous $g$
\eqn{Leb}
$$
V(g_{j},[0,1]^2)\leq A(g, [0,1]^2)\leq
V(g_{1},[0,1]^2)+V(g_{2},[0,1]^2)+V(g_{3},[0,1]^2).
$$

Let us note that these results are  highly nontrivial.  For example it is
possible to construct continuous $g$ such that
$A(g,[0,1]^2)$ is much smaller than  $\haus^2(g([0,1]^2))$ (which may be even
infinite) but the result \eqref{Leb} is still true. Further for the validity we
need only continuity of $g$ and we do not need to assume that $A(g,
[0,1]^2)<\infty$. However, this is only  known to hold for two dimensional
surfaces in
$\er^3$ and for higher dimensions the assumption about the finiteness of the
Lebesgue area might be needed.

\subsection{$BV$ functions and the coarea formula}
Let $\Omega\subset\rn$ be an open domain.
A function $h\in L^1(\Omega)$ is of bounded variation, $h\in BV(\Omega),$
if the distributional partial derivatives of $h$ are measures with finite total
 variation in $\Omega$, i.e. there are Radon (signed) measures $\mu_1,\ldots,\mu_n$ defined in $\Omega$
so that for $i=1,\ldots,n,$ $|\mu_i|(\Omega)<\infty$ and
$$\int_{\Omega} hD_i\varphi\; dx=-\int_{\Omega}\varphi\; d\mu_i$$
for all $\varphi \in C^{\infty}_0(\Omega).$
We say that $f\in L^1(\Omega,\er^n)$ belongs to $BV(\Omega,\er^n)$ if the
coordinate functions of $f$
belong to $BV(\Omega)$.

Let $\Omega\subset\rn$ be open set and $E\subset\Omega$ be measurable. The
perimeter of $E$ in
$\Omega$ is defined as a total variation of $\chi_E$ in $\Omega$, i.e.
$$
P(E,\Omega)\colonequals  \sup\Bigl\{\int_{E}\operatorname{div}\varphi\; dx:\
\varphi\in C_0^1(\Omega),\ \|\varphi\|_{\infty}\leq 1\Bigr\}.
$$

We will need the following coarea formula to
characterize $BV$ functions (see \cite[Theorem 3.40]{AFP}).

\prt{Theorem}
\begin{proclaim}\label{thm:CoareaFormulaForBV}
  Let $\Omega\subset\rn$ be open and $u\in L^1(\Omega)$. Then we have
  \eqn{coarea}
  $$
  |Du|(\Omega)=\int_{-\infty}^{\infty}P(\{x\in \Omega:\ u(x)>t\},\Omega)\; dt.
  $$
  In particular, $u\in BV(\Omega)$ if and only if the integral on the righthand side is finite.
\end{proclaim}

It is well-known (see e.g.\  \cite[Proposition 3.62]{AFP}) that for the coordinate
functions of a homeomorphism $f \colon \Omega\to\rn$ we have
\eqn{perimeter}
$$
P\bigl(\{x\in \Omega:\ f_i(x)>t\},\Omega\bigr)\leq \haus^{n-1}\bigl(\{x\in \Omega:\ f_i(x)=t\}\bigr).
$$

Moreover, we have the following  version of coarea formula for continuous $BV$
functions by Federer \cite[Theorem 4.5.9 (13) and (14) for $k\equiv 1$]{Fe}.
\prt{Theorem}
\begin{proclaim}\label{coarea2}
  Let $\Omega\subset\rn$ be open and $u\in BV(\Omega)$ be continuous. Then we have
  % \eqn{coarea}
  $$
  |Du|(\Omega)=\int_{-\infty}^{\infty}\haus^{n-1}(\{x\in \Omega:\ u(x)=t\})\; dt.
  $$
\end{proclaim}

\subsection{BVL condition}\label{sec:BVL}
Let $\Omega \subset\rn$ be open and $f\in L^1(\Omega)$. It is well-known that $f\in BV(\Omega)$ if and only if it satisfies the BVL condition,
i.e.\ it has bounded variation on $\mathcal{L}^{n-1}$ a.e.\ line parallel to the coordinate axes, and the variation along
these lines is integrable (see e.g.\ \cite[Remark 3.104]{AFP}). As a corollary we obtain that a $BV$ function of
$n$-variables is a $BV$ function of $(n-1)$-variables on $\mathcal{L}^1$ a.e.\ hyperplane parallel to coordinate axis.

For example for $n=2$ and $f\in BV((0,1)^2)$ we have that the function
$$
f_x(y)\colonequals f(x,y)
$$
has bounded (one-dimensional) variation for a.e.\ $x\in(0,1)$. Moreover,
\eqn{eq:y-DirectionVariation}
$$
\int_0^1 |D f_x((0,1))|\; dx=|D_2 f|((0,1)^2)
$$
where $|D f_x|$ denotes the (one-dimensional)  total variation of $f_x$ and
$|D_2 f|$ denotes the total variation of the measure
$\frac{\partial{f}}{\partial y}$. Similar identity holds for $f_y(x)\colonequals
f(x,y)$ and $D_1 f$.

\subsection{Convergence of $BV$ functions}
\label{sec:Convergence}

In dimension two, the boundary of a ball $B(x,r)$
is a curve and we will tacitly assume that it is always
parametrized with the path
\begin{align*}
  \beta \colon [0,2\pi] \to \er^2, \quad
  \beta(t)
  = (x_1 + r \cos t , x_2 + r \sin t).
\end{align*}
Thus when we speak of the \emph{length} $\ell(\partial B(x,r))$
of the boundary of a ball or its image $f (\partial B(x,r))$
under a mapping $f$, we mean the length of the curve $\beta$
or $f \circ \beta$, respectively. Note that the length of a
path $\gamma \colon [0,1] \to \er^2$ equals
\begin{align*}
  \ell(\gamma)
  \colonequals
  \left\{
  \sum_{j=1}^{k} d(\gamma(t_{j-1}),\gamma(t_j)) :
  0 = t_0 \leq \ldots \leq t_k = 1
  \right\},
\end{align*}
from which we immediately see that if $f_j \to f$ uniformly, then
\begin{align*}
  \lim_{j\to\infty} \ell(f_j \circ \beta) \geq \ell(f \circ \beta).
\end{align*}
Similarly we also assume line segments in the plane to be equipped with
a path parametrization and to have similar length convergence properties.

By the results in the previous subsection \ref{sec:BVL} we know that the restriction of a $BV$ mapping
$f \in BV(\er^2,\er^2)$
to $\mathcal{L}^1$ a.e.\ line segment in the plane
is again a $BV$ mapping;
i.e.\ for all $a>b$ and a.e.\ $x \in \er$, the restriction
$f_x \colonequals  f|_{\{ x \} \times (a,b)}$ is a $BV$ mapping and
\begin{align*}
  \mathcal{H}^1(f (I_x))
  \leq \ell(f (I_x))
  = | Df_x |(I_x)
  < \infty,
\end{align*}
where $| Df_x |$ denotes the one-dimensional total variation of $f_x$.

A similar result holds also for $\mathcal{H}^1$ a.e.\ radius of a sphere: given a point $x$, the restriction
of $f$ to $\partial B(x,r)$ is $BV$ for $\mathcal{H}^1$ a.e.\ radius $r > 0$.
This especially implies that for such radii,
\begin{align*}
  \mathcal{H}^1(f (\partial B(x,r)))
  \leq \ell(f (\partial B(x,r)))
  = | D(f_r) |(\partial B(x,r))
  < \infty,
\end{align*}
where $f_r \colonequals  f|_{\partial B(x,r)}$ and $|D(f_r)|$ denotes the one-dimensional total variation of $f_r$.
Furthermore, similarly as in \eqref{eq:y-DirectionVariation}, we have
\eqn{eq:TangentialVariation}
$$
  \int_0^r \mathcal{H}^1(f (\partial B(x,s))) \; ds
  \leq \int_0^r| D(f_s) |(\partial B(x,s)) \;ds
  \leq |Df|\left(B(x,r)\right).
$$

Recall the weak* convergence of $BV$ mappings.
\prt{Definition}
\begin{proclaim}
  We say that a sequence $(f_j)$ of $BV$ mappings
  \emph{weakly* converges to $f$ in $BV$}, if $f_j \to f$ in $L^1$
  and $Df_j$ weakly* converge to $Df$, i.e.\
  \begin{align*}
    \lim_{j \to \infty} \int_\Omega \phi \; dDf_j
    = \int_\Omega \phi \; dDf
  \end{align*}
  for all $\phi \in C_0(\Omega)$.
\end{proclaim}

The following result from \cite[p.125, Proposition 3.13.]{AFP}
gives a characterization for weak* convergence in $BV$. Note especially that since  the
mollifications of continuous $BV$ functions converge uniformly,
they especially converge locally in $L^1$, so in this case the
boundedness of the sequence in $BV$-norm gives weak* convergence
for the derivatives.
\prt{Proposition}
\begin{proclaim}
  \label{prop:BV-Convergence}
  Let $f_j$ be a sequence of $BV$ mappings $\Omega \to \er^2$. Then
  $f_j$ weakly* converges to a $BV$ mapping $f \colon \Omega \to \er^2$
  if and only $f_j \to f$ in $L^1$ and $\sup |Df_j|(\Omega)<\infty.$
\end{proclaim}

\section{Properties of $BV$ mappings}

In the proof of Theorem \ref{regularity} we use some ideas of \cite[proof of Theorem 1.7]{DS}.
In particular we use the following observation based on the coarea formula (Theorem \ref{thm:CoareaFormulaForBV}).

\prt{Theorem}
\begin{proclaim}\label{abstract}
  Let $\Omega\subset\er^n$ be a domain and $f\in BV_{\loc}(\Omega,\er^n)$ be a homeomorphism.
  Then the following measure on $\Omega$ is finite
  $$
  \mu(A)=\sum_{i=1}^n\int_{-\infty}^{\infty}\haus^{n-1}\bigl(f(\{x\in A:\ x_i=t\})\bigr)\; dt
  $$
  if and only if $f^{-1}\in BV_{\loc}(f(\Omega),\er^3)$.
  In addition, $f(\mu)$ is absolutely continuous with respect to the Lebesgue
  measure if and only if $f^{-1}\in W^{1,1}(f(\Omega),\er^3).$
\end{proclaim}
\begin{proof}
  Assume that $\mu$ is a finite measure. By Theorem \ref{thm:CoareaFormulaForBV} and the perimeter inequality \eqref{perimeter} we have
  \eqn{www}
  $$
  \begin{aligned}
    |Df^{-1}|(f(\Omega))&\approx\sum_{i=1}^n |(Df^{-1})_i|(f(\Omega))\\
    &= C\sum_{i=1}^n\int_{-\infty}^{\infty}P(\{y\in f(\Omega):\ (f^{-1})_i(y)>t\},f(\Omega))\; dt\\
    &\leq C\sum_{i=1}^n\int_{-\infty}^{\infty}\haus^{n-1}\bigl(\{y\in f(\Omega):\ (f^{-1})_i(y)=t\}\bigr)\; dt\\
    &=C\sum_{i=1}^n\int_{-\infty}^{\infty}\haus^{n-1}\bigl(f(\{x\in\Omega:\ x_i=t\})\bigr)\; dt<\infty.\\
  \end{aligned}
  $$
  and thus $f^{-1}\in BV_{\loc}$.

  If $f^{-1}\in BV_{\loc}$, then by Theorem \ref{coarea2} we know that the only
  inequality in the above computation \eqref{www} is actually equality  and we get
  $\mu\in \M(\Omega)$.

  Let us consider now the final claim. We have to show that
$\abs{Df^{-1}}(E)<\epsilon$
  if $\abs{E}<\delta.$
  Given $\epsilon$ we choose $\delta>0$ from the absolute continuity of measure $f(\mu)$.
  By approximation we may assume that $E$ is open and $|E|<\delta$. The
  definition of $\mu,$ assumed absolute continuity of $f(\mu)$ and \eqref{www} (with $E$ instead of $\Omega$) imply
  $$
  \abs{Df^{-1}}(E)\leq \mu(f^{-1}(E))<\epsilon.
  $$
  If we know that $f^{-1}\in W^{1,1}$ then we have only equalities in \eqref{www} and we easily obtain that $f(\mu)$ is absolutely continuous with respect to the  Lebesgue measure.
\end{proof}

We next show that for a mollification of a continuous $BV$ mapping,
the convergence is inherited to a.e.\ circle in a weak sense.
\prt{Proposition}
\begin{proclaim}\label{prop:BoundaryLimits}
  Let $f \colon \Omega \subset \er^2 \to \er^2$ be
  a continuous $BV$ mapping and $(\rho_k)$ a sequence of
  mollifiers. Denote $f^k \colonequals  f \ast \rho_k$,
  $f^k_{c,s} \colonequals  f^k|_{\partial B(c,s)}$ and
  $f_{c,s} \colonequals  f|_{\partial B(c,s)}$.
  Then for any point $z \in \er^2$ we have
  $$
  \lim_{k \to \infty}| D f^k_{z,r}|(\partial B(z,r))
  = |Df_{z,r}|(\partial B(z,r))
  < \infty
  $$
  and
  $$
  D (f^k_{z,r}) \overset{w*}{\rightharpoonup} D (f_{z,r}),
  $$
  for $\mathcal{H}^1$ a.e.\ radius $r > 0$
  such that $\overline{B(z,r)} \subset \Omega$.
\end{proclaim}

\begin{proof}
  Since the claim is local  it suffices, after a smooth
  change of local coordinates, to show that for a continuous
  $BV$ mapping $f \colon (0,1)^2 \to \er^2$ we have
  \eqn{eq:BoundaryIntegralLimit}
  $$
  \lim_{k \to \infty} | D f_x^k|(I_x)
  = | D f_x|(I_x)
  < \infty,
  $$
  and
  \eqn{eq:BoundaryWeakConvergence}
  $$
  D (f^k_x) \overset{w*}{\rightharpoonup} D (f_x),
  $$
  on $\mathcal{H}^1$ almost every line segment
  $I_x \colonequals  \{ x \} \times (0,1)$,
  where $f_x^k \colonequals  f^k|_{I_x}$ and
  $f_x \colonequals  f|_{I_x}$.

  We start by proving \eqref{eq:BoundaryIntegralLimit}.
  By the results in Section \ref{sec:Convergence},
  for $\mathcal{H}^1$ a.e.\ $x \in (0,1)$
  we have
  $$
  |Df_x^k|(I_x)
  = \ell(f^k(I_x))
  < \infty
  \quad \text{ and } \quad
  |Df_x|(I_x)
  = \ell(f(I_x))
  < \infty,
  $$
  so since $f_x^k \to f_x$ uniformly we see by the notions of Section \ref{sec:Convergence} that
  \begin{align*}
    \lim_{k \to \infty} |Df_x^k|(I_x)
    =\lim_{k \to \infty} \ell(f^k(I_x))
    \geq \ell(f(I_x))
    = |Df_x|(I_x).
  \end{align*}
  Thus to prove \eqref{eq:BoundaryIntegralLimit} it suffices to show that
  for a.e.\ $x \in (0,1)$,
  \begin{align*}
    \lim_{k \to \infty} |D f_x^k|(I_x)
    \leq |Df_x|(I_x).
  \end{align*}

  Suppose this is not true, whence there exists
  $\delta > 0$ such that the set
  $$
  J\colonequals \Bigl\{ x\in (0,1):\ \lim_{k \to \infty} |Df_x^k|(I_x)
  > (1+\delta)|D f_x|(I_x) \Bigr\}
  $$
  has positive 1-measure.
  Fix a Lebesgue point $x_0 \in (0,1)$ of $J$.
  By the Lebesgue density theorem we may assume $x_0$ to be such that
  $$
  \lim_{r\to 0}\frac{1}{2r}\int_{x_0-r}^{x_0+r}\bigl||Df_x|(I_x)-|Df_{x_0}|(I_{x_0}) \bigr|
  =0.
  $$
  Choose $\eta>0$ such that
  \eqn{triv}
  $$
  (1+\delta)\frac{1-\eta}{2}\bigl((1-\eta)2-\eta\bigr)>1.
  $$
  Fix $r>0$ for which
  \eqn{lot}
  $$
  \begin{aligned}
    (i)\ &|D f|\bigl(\partial ((x_0-r,x_0+r) \times (0,1))\bigr) = 0,\\
    (ii)\ &|J\cap (x_0-r,x_0+r)|\geq (1-\eta) 2r,\\
    (iii)\ &\int_{x_0-r}^{x_0+r}\bigl||Df_x|(I_x)-|Df_{x_0}|(I_{x_0}) \bigr| \;
dx< \eta |D f_{x_0}|(I_{x_0}) r, \textrm{ and }\\
    (iv)\ &|D f_{x_0}|(I_{x_0})>\frac{1-\eta}{2r} \int_{x_0-r}^{x_0+r}|Df_x|(I_x)\; dx
    .\\
  \end{aligned}
  $$

  As remarked in Section \ref{sec:BVL}, $D_2f$ is a finite Radon measure. Thus
applying
  \cite[Proposition 2.2.(b), p.42]{AFP} for the mollification of its total variation $|D_2 f|$ and
  using the fact that the measure is Borel regular, we see that
  \begin{align*}
    \lim_{k \to \infty} | D_2 f^k |(U)
    \leq \lim_{k \to \infty} | D_2 f |\left(U+B(0,k \inv)\right)
    = |D_2 f|(\overline{U}).
  \end{align*}
  for any Borel set $U$.
  By setting $U \colonequals  (x_0-r,x_0+r) \times (0,1)$
  we have by $(i)$ that $|D f| (\partial U) = 0$,
  and so also $|D_2 f| (\partial U) \leq |D f|(\partial U) = 0$. Thus
  \begin{align}\label{eq:VariationLimit}
    \lim_{k \to \infty} |D_2 f^k|(U)
    \leq |D_2 f|(\overline{U})
    = |D_2 f|(U).
  \end{align}
  On the other hand by using Fatou's lemma, the definition of $J$, \eqref{lot} $(iii)$, $(ii)$, $(iv)$ and \eqref{triv},
  \begin{align*}
    \lim_{k \to \infty} |D_2 f^k|(U)
    & \geq \lim_{k \to \infty} \int_{(x_0-r,x_0+r)\cap J} |D f^k_x|(I_x) \; dx  \\
    & \geq (1+\delta) \int_{(x_0-r,x_0+r)\cap J} |D f_x|(I_x) \; dx  \\
    &\geq (1+\delta) \Bigl[\int_{(x_0-r,x_0+r)\cap J} |D f_{x_0}|(I_{x_0})\; dx-\eta |D f_{x_0}|(I_{x_0}) r\Bigr]\\
    &\geq (1+\delta)|D f_{x_0}|(I_{x_0})\bigl((1-\eta)2r-\eta r\bigr)\\
    &\geq (1+\delta) \frac{1-\eta}{2r} \int_{x_0-r}^{x_0+r}|Df_x|(I_x)\; dx\bigl((1-\eta)2r-\eta r\bigr)\\
    &> |D_2 f^k|(U).
  \end{align*}
  This contradicts \eqref{eq:VariationLimit} and so \eqref{eq:BoundaryIntegralLimit} holds.

  To prove \eqref{eq:BoundaryWeakConvergence} we note that for $\mathcal{H}^1$ a.e.\
  line segment $I_x$ the $BV$ mappings $f_x^k \colon I_x \to \er^2$
  converge uniformly to the continuous $BV$ mapping $f_x \colon I_x \to \er^2$.
  Furthermore they form a bounded sequence with respect to the $BV$ norm, and
  thus by Proposition \ref{prop:BV-Convergence} they converge weak* in $BV$. This implies
  \eqref{eq:BoundaryWeakConvergence} and the proof is complete.
\end{proof}

\section{Degree Theorem for continuous $BV$ planar mappings}

The aim of this section is to prove the following analogy of the
change of variables formula for the distributional Jacobian
in two dimensions.
A similar statement was shown before in \cite{CL} for mappings that satisfy $J_f>0$ a.e.\ and that are one-to-one and in \cite{DHMS} for open and discrete mappings. Here we generalize this result to mappings where the Jacobian can change the sign but we restrict our attention to planar mappings only.

\prt{Theorem}
\begin{proclaim}
  \label{thm:degree}
  Let $f \colon \er^2 \to \er^2$ be a continuous
  $BV$ mapping such that the distributional
  Jacobian $\J f$ is a signed Radon measure.
  Then for every $ x \in \er^2$ we have
  \begin{align}\label{eq:DegreeFormula}
    \int_{\er^2} \deg(f,B(x,r),y) \; dy
    = \J_f \left( B(x,r)\right)
  \end{align}
  for a.e.\ $r > 0$.
\end{proclaim}

Before the proof Theorem \ref{thm:degree} we prove the following important corollary,
which is one of the main tools in the proof of Theorem \ref{regularity}.

\prt{Proposition}
\begin{proclaim}
  \label{coro:degreeWithAbsoluteValues}
  Let $f \colon \er^2 \to \er^2$ be a continuous
  $BV$ mapping such that the distributional
  Jacobian $\J f$ is a signed Radon measure and such that $V(f,\er^2)<\infty$.
  Then for every $ x \in \er^2$ we have
  \begin{align}\label{eq:DegreeCorollary}
    \int_{\er^2} \left| \deg(f,B(x,r),y) \right| \; dy
    \leq
    \left|\J f\right| \left( B(x,r)\right)
  \end{align}
  for a.e.\ $r > 0$.
\end{proclaim}
\begin{proof}

Let us note that the previous theorem holds not only for balls but also for a.e. cube $Q(x,r)$. From the previous theorem we know that the set 
$$
\bigl\{[x,r]:\ Q(x,r)\text{ is good for }\eqref{eq:DegreeFormula},\ \mathcal{L}^2(f(\partial Q(x,r)))=0\text{ and }|\J_f|(\partial Q(x,r))=0\bigr\}. 
$$
has full $\mathcal{L}^{n+1}$ measure. 
It follows that for a.e. $r>0$ we have that $Q(x,r)$ is good for a.e. $x\in \Omega$ with $r<\dist(x,\partial \Omega)$. 
Hence we can fix $r_0>0$ such that all $r_k=r_0 2^k$, $k\in\zet$, are good for every 
$x\in\Omega\setminus N_0$ with $|N_0|=0$. Hence we can fix $x_0\in\rn$ and a dyadic grid 
\eqn{defg0}
$$
G_0:=\bigl\{x_0+2^k Q(y_i,r_0),\ y_i\in\zet^n\bigr\}, 
$$
such that for all cubes from the grid inside $\Omega$ we have 
\eqn{key}
$$
\int_{\er^2} \deg(f,Q,y) \; dy= \J_f \left( Q\right). 
$$  
 and 
$\mathcal{L}^2(f(\partial Q))=|\J_f|(\partial Q)=0$ for every $Q\in G_0$, $Q\subset \Omega$. It is enough to choose any 
$$
x_0\notin \bigcup_{k=1}^{\infty}\Bigl(\bigcup_{y_i\in\zet^n} 2^k \{-y_i\}+N_0\Bigr). 
$$

Let us fix a cube $Q\subset\Omega$. Instead of proving \eqref{eq:DegreeCorollary} for a ball we prove it for $Q$ which is equivalent. 
Analogously to \eqref{recalldef} we define (see \cite[9.10]{C})
$$
U(f,Q)\colonequals  \sup_S\left\{\sum_{\pi\in S}\Bigl|\int_{\er^2}
\deg(f,\pi,y)\; dy\Bigr|\right\},
$$
where $S$ is any finite system of nonoverlapping simple open polygonal regions in
$Q$. 
By \cite[12.9 Theorem (iii)]{C} we know that
$$
U(f,Q)=V(f,Q).  
$$
Note that for the validity of this identity we need the additional assumption $V(f,Q)<\infty$ as it is the assumption of \cite[12.9 Theorem (iii)]{C} (it is stated in \cite{C} as $f$ is plane $BV$ but in the notation of the book it means exactly $V(f,Q)<\infty$).
By \cite[12.9 Theorem (i) and 12.6]{C} we know that there is a sequence of figures $F_n$ in $Q$ that consists of disjoint cubes from our dyadic grid $Q_{i,n}$, $F_n=\bigcup_i Q_{i,n}$, so that 
$$
U(f,Q)=\lim_{n\to\infty} \sum_{i}\Bigl|\int_{\er^2}\deg(f,Q_{i,n},y)\; dy\Bigr|. 
$$

Now we can easily estimate with the help of \eqref{key} and definition of total variation
$$
\begin{aligned}
\int_{\er^2} \left| \deg(f,Q,y) \right| \; dy&\leq V(f,Q)=U(f,Q)\\
&=\lim_{n\to\infty} \sum_{i}\Bigl|\int_{\er^2}\deg(f,Q_{i,n},y)\; dy\Bigr|\\
&\leq\lim_{n\to\infty} \sum_{i}\Bigl|\djac(Q_{i,n})|\\
    &\leq \left|\J f\right| \left( Q\right).\\
\end{aligned}    
$$
\end{proof}

The proof of Theorem \ref{thm:degree} requires several auxiliary results.
We begin with the following degree convergence lemma; compare to Lemma
\ref{lemma:TopologicalDegree}.
\prt{Lemma}
\begin{proclaim}\label{lemma:DegreeLimit}
  Let $f \colon \Omega \subset \er^2 \to \er^2$ be
  a continuous $BV$ mapping and let $(f^k)$ be
  mollifications of $f$. Then for any point $x \in \er^2$ and a.e.\ radius $r > 0$
  we have
  \begin{align}\label{eq:DegreeLimit}
    \lim_{k \to \infty} \int_{\er^2} \deg(f^k,B(x,r),y) \; dy
    = \int_{\er^2} \deg(f,B(x,r),y) \; dy.
  \end{align}
\end{proclaim}
\begin{proof}
  Let $x_0 \in \er^2$.
  By the BVL properties remarked in Section \ref{sec:BVL} we know that
  for almost every radius $r>0$ the length of $f(\partial B(x_0,r))$
is finite, i.e.
  %$\ell(f \partial B(x_0,r))< \infty$.
  %\lim_{k \to \infty}| D f^k_{z,r}|(\partial B(z,r))
  $|Df_{z,r}|(\partial B(z,r))<\infty$ and that the claim of Proposition
\ref{prop:BoundaryLimits} holds.
  Fix such a $r_0$, and set $B_0 \colonequals  B(x_0,r_0)$.

  We first define
  $$
  F^k(y)= \deg(f^k,B_0,y)
  \quad
  \text{ and }
  \quad
  F(y)= \deg(f,B_0,y),
  $$
  whence
  $$
  \int_{\er^2} |\deg(f^k,B_0,y)| dy=\norm{F^k}_1
  \quad
  \text{ and }
  \quad
  \int_{\er^2} |\deg(f,B_0,y)| dy=\norm{F}_1
  $$
  and we have to show that $F^k\rightarrow F$ in $L^1.$
  To show this we use compactness of $BV$.
  First we show that $F^k$ is bounded sequence in $BV$-norm.

  It is easy to see that the variation measure of $F^k$ is supported only on the
  curve $f^k(\partial B_0)$. Furthermore, since
  $f^k(\partial B_0)$ is rectifiable, $\mathcal{H}^1$-a.e.\ point
  is on the boundary of at most two components of $\er^2 \setminus f^k(\partial B_0)$.
  In such a situation, if the value of $F^k$ differs by $N$ on these two components,
  the image $f^k(\partial B_0)$ must cover this joint boundary at least
  $N$ times. Thus the total variation of $F^k$ is in
  fact bounded by
  $$\ell\bigl(f^k(\partial B_0)\bigr)=\abs{Df^k_{x_0,r_0}}(\partial
  B_0).$$
  Since the radius $r_0$ was chosen such that
  Proposition \ref{prop:BoundaryLimits} holds, we have
  $$
  \abs{Df^k_{x_0,r_0}}(\partial B_0)\rightarrow \abs{Df_{x_0,r_0}}(\partial B_0)$$
  and so $\abs{DF^k}$ is
  uniformly bounded.
  Furthermore the boundedness of the sequence $(F^k)$ in $L^1$ follows from the Sobolev inequality \cite[Theorem
  3.47]{AFP}.
  Thus, the compactness theorem in \cite[Theorem 3.23]{AFP} implies that there exists a
  subsequence $(F^{k(j)})$ which converges in $L^1$ to a function $G$.

  We will show that $G=F,$ which implies that the original sequence $F^k$
  converges to $F$ in $L^1,$ as every converging subsequence must converge to $F.$
  Assume that $F\neq G$ on a set $A$ with positive Lebesgue measure. Since
  $f(\partial B_0)$ has finite 1-Hausdorff measure we find with the Lebesgue
  density theorem   $z\in \er^2\setminus f(\partial B_0),$ which is a
  density point  of $A$ with $G(z)\neq F(z).$
  For some very small ball $B_z$ centered at $z$ we have
  $$
  \abs{\int_{B_z}G-F}>0
  $$
  and $B_z$ is compactly contained in some component of $\er^2\setminus f(\partial
  B_0).$
  Now recall that  $f^k$ converge uniformly to $f. $
  When $\norm{f^k-f}<\dist(B_z,f(\partial B_0))$ we have by basic
  properties of the degree (see \cite[Theorem 2.3.]{FG})
  $$
  F^k(y)=\deg(f^k,B_0,y)=\deg(f,B_0,y)=F(y)
  $$
  for every point  $y\in B_z.$ This is a contradiction with $L^1$ convergence and
  the definition of $B_z$. Thus the original claim follows.
  \end{proof}

The proof of the previous lemma goes through also with
absolute values of the degrees. We record this observation
as the following corollary even though we will not be using
it in this paper.

\prt{Corollary}
\begin{proclaim}
  \label{corollary:DegreeLimit}
  Let $f \colon \Omega \subset \er^2 \to \er^2$ be
  a continuous $BV$ mapping and $(f^k)$ a mollification of $f$.
  Then for any point $x \in \er^2$ and a.e.\ radius $r > 0$
  we have
  \begin{align}\label{eq:DegreeLimit2}
    \lim_{k \to \infty} \int_{\er^2} |\deg(f^k,B(x,r),y)| \; dy
    = \int_{\er^2} |\deg(f,B(x,r),y)| \; dy.
  \end{align}
\end{proclaim}

The following Proposition \ref{prop:degree} is essentially
a $BV$-version of \cite[Proposition 2.10]{HMC}.
For smooth mappings the identity \eqref{eq:DegreeEquation}
follows in a more general form with smooth test functions $g \in C^\infty(\Omega,\er^2)$ by combining
the Gauss-Green theorem and the area formula in a ball $B$:
\begin{align}
\label{MSTformula}
  \int_{\partial B} \left\langle ( g(f(x)) \cdot \cof D f(x) , \nu  \right\rangle \; d\mathcal{H}^1(x)
  &= \int_{B} \operatorname{div} g (f(y)) J_f(y)\; dy \\
  &= \int_{\er^2}\operatorname{div} g(y)\deg(f, B, y)\; dy\nonumber,
\end{align}
where $\nu$ denotes the unit exterior normal to $B$ and $\cof D f(x)$ denotes the cofactor matrix, i.e. the matrix of $(n-1)\times(n-1)$ subdeterminants with correct signs.
For more details for the general setting we refer to M\"uller, Spector and Tang; in \cite[Proposition 2.1]{MST}
they prove the claim for continuous $f \in W^{1,p}$, $p > n-1$ and $g \in C^1$.
We need the identity only in the case of $g(x_1,x_2)=[x_1,0]$. In this case the
integrand on left hand side of \eqref{MSTformula} reduces to
\begin{equation}
 f_1\left\langle Df_2, \nu_t \right\rangle,
\end{equation}
where $\nu_t$ is the unit tangent vector of $\partial B.$ Thus in the $BV$
setting it is natural to replace the left hand side of \eqref{MSTformula}
with
$$
\int_{\partial B} f_1 d(Df|_{\partial B}).
$$
since by Section \ref{sec:BVL} $f$ is one dimensional $BV$-function on almost
every sphere centered at any given point.

\prt{Proposition}
\begin{proclaim}\label{prop:degree}
  Let $\Omega\subset\er^2$ be a domain and let $f \colon \Omega \to \er^2$ be a continuous $BV$ mapping.
  Then for every $c\in\er^2$ and a.e.\ $r>0$ such that $B\colonequals B(c,r)\subset\Omega$ we have
  \begin{align}\label{eq:DegreeEquation}
\int_{\partial B}f_1 d(Df|_{\partial B})
    =\int_{\er^2}\deg(f, B, y)\; dy.
  \end{align}

\end{proclaim}

\begin{proof}
  We prove the claim by approximating $f$ with a sequence of mollifiers $(f^k)$,
  showing that $[f^k_1,0] \cdot \cof Df^k$ convergences weakly* to $f_1
d(Df|_{\partial B})$ and combining
  this with Lemma \ref{lemma:DegreeLimit}.

  Let us fix $r>0$ such that $B(c,r)\subset\subset\Omega$ and the conclusion of Lemma \ref{lemma:DegreeLimit} and Proposition \ref{prop:BoundaryLimits} hold for this radius.
  Now for every $k$ with $B(c,r+\frac{1}{k})\subset\subset\Omega$ we have
  $f^k\in C^\infty(B,\er^2)$.
  Since $f$ is continuous, $f^k \to f$ uniformly.
  Clearly
$$
   \int_{\partial B} \bigl( f^k_1 Df_2^k|_{\partial B} - f_1 Df|_{\partial B} \bigr)
    =     \int_{\partial B}
    \left(f_1 Df_2^k|_{\partial B} - f_1  Df|_{\partial B} \right)
    -\int_{\partial B}
    \bigl(f_1^k - f_1 \bigr)Df_2^k|_{\partial B}.
$$
  We next note that by Proposition \ref{prop:BoundaryLimits}, $Df^k|_{\partial B} \to Df|_{\partial B}$ with respect to the weak* convergence
  and $\|f_1^k - f_1\|_\infty \to 0$ by the uniform convergence of $(f^k)$.
  Thus both terms of the right hand side converge to
  zero as $k \to \infty$. It follows that
  \eqn{eq:CofactorLimit}
  $$
  \begin{aligned}
    \lim_{k \to \infty}\int_{\partial B}  &\bigl\langle [f_1^k(x),0] \cdot \cof Df^k(x), \nu \bigr\rangle\, d \mathcal H^{1} (x)=\\
    &=\lim_{k \to \infty}\int_{\partial B}  f_1^k(x) Df^k|_{\partial B}(x) d
\mathcal H^{1} (x)
    =\int_{\partial B} f_1(x) d(Df|_{\partial B}(x)).
    \end{aligned}
  $$

  On the other hand, since the mappings $f^k$ are smooth we have by e.g.\ \cite[Proposition 2.1]{MST}
  that
  \begin{align*}
    \int_{\partial B}  \bigl\langle [f_1^k(x),0] \cdot \cof Df^k(x), \nu \bigr\rangle\; d \mathcal H^{1} (x)
    =\int_{\er^2}\deg(f^k, B, y) dy.
  \end{align*}
  Combining this with \eqref{eq:CofactorLimit} and Lemma \ref{lemma:DegreeLimit} gives the claim.
\end{proof}

We are now ready to prove the main result of this section, Theorem \ref{thm:degree}. In its proof we use some ideas from \cite{CL} and \cite{M}.

\begin{proof}[Proof of Theorem \ref{thm:degree}.]
  We recall the definition of distributional Jacobian for any $\varphi \in C^{\infty}_0(\Omega)$
  \eqn{rrr}
  $$
  %\begin{aligned}
    \mathcal J_f(\varphi)
    =-\int_{\Omega}f_1(x) J(\varphi(x),f_2(x))\; dx
    =\int_{\Omega} \bigl\langle [f_1(x),0] \cdot \cof Df(x),D\varphi(x) \bigr\rangle\; dx.
  %\end{aligned}
  $$

  Let us pick a ball $B\colonequals B(y,r)\subset\Omega$ such that $|Df|(\partial B)=0$.
  Furthermore, by the Lebesgue theorem we may assume that
  \eqn{lebpoint}
  $$
  \lim_{\delta\to 0}\frac{1}{\delta}\int_{r -\delta}^r\Bigl||Df_{y,s}|(\partial
B(y,s))-|Df_{y,r}|(\partial B(y,r))\Bigr|\; ds=0
  $$
  where $f_{y,s}\colonequals f|_{\partial  B(y,s)}$ and $|Df_{y,s}|$ is the
corresponding (one-dimensional) total variation.
  Let us fix $\psi\in C^\infty(\er,[0,1])$  such that $\psi(s)\equiv 1$ for $s<0$
and $\psi(s)\equiv 0$ for $s>1$. For $0<\delta< r$ we set
  $$
  \Phi_\delta(s)=\psi\Bigl(\frac{s-(r-\delta)}{\delta}\Bigr), \text{ i.e. }
  \Phi_{\delta}(s)=\begin{cases}
    1\text{ for }s\leq r-\delta.\\
    0\text{ for }s\geq r.\\
  \end{cases}\text{ and }|\Phi'_{\delta}|\leq \frac{C}{\delta}.
  $$
  As the distributional Jacobian is a Radon measure and $|Df|(\partial B)=0$ we obtain
  \begin{align}\label{eq:InnerRegularity}
    \mathcal J_f(B(y,r))
    = \lim_{\delta\to 0+} \int_\Omega \Phi_\delta(|x-y|)\; d\mathcal J_f(x).
  \end{align}
  By \eqref{rrr} for $\varphi=\Phi_{\delta}({ |x-y|})$ and Proposition \ref{prop:degree} we have
    \eqn{ahoj}                                 %
  $$
  \begin{aligned}
    \int_{\Omega}  \Phi_\delta(|x-y|)\; d\mathcal J_f(x)
    &= \int_{\Omega} \bigl\langle [f_1(x),0] \cdot \cof Df(x), D\Phi_{\delta}(|x-y|)   \bigr\rangle\, dx\\
    &= \int_{r-\delta}^{r}  \int_{\partial B(y,s)} f_1(x) \Phi'_{\delta}(s) d(Df|_{\partial B}(x))\\
&= \int_{r-\delta}^{r} \Phi'_{\delta}(s) \int_{\er^2} \deg (f, B(y,s), z)\;  dz\; ds.\\
  \end{aligned}
  $$

  We next show that the integral on the right hand side of \eqref{ahoj} converges as $\delta \to 0$.
  For all $y \in \Omega$ and $\delta > 0$ small enough we set
  \begin{align*}
    f_s(x)
    =f\left(\frac{s}{r}(x-y)+y\right).
  \end{align*}
  Note that with this notation
  \eqn{444}
  $$
  \begin{aligned}
    \int_{r-\delta}^{r} \Phi'_{\delta}(s) \int_{\er^2} \deg (f, B(y,s), z)\;  dz\; ds
    &= \int_{r-\delta}^{r} \Phi'_{\delta}(s) \int_{\er^2} \deg (f_s, B(y,r), z)\;  dz\; ds
  \end{aligned}
  $$
  and the right hand side is a type of average integral as $\int_{r-\delta}^r \Phi'_\delta=1$.

We have a fixed mapping $f|_{\overline{B(y,r)}}$ with $|Df_{y,r}|(\partial B(y,r))<\infty$.
We claim that
given $\epsilon>0$ we can find $\eta>0$ such that for every continuous mapping
$g|_{\overline{B(y,r)}}$ we have
\eqn{reallyclose}
$$
\begin{aligned}
\|f-g\|_{L^{\infty}(\partial B)}&<\eta\text{ and }\bigl||Df_{y,r}|(\partial B(y,r))- |Dg_{y,r}|(\partial B(y,r))\bigr|<\eta\Rightarrow\\
&\Rightarrow \Bigl|\int_{\er^2} \deg(f,B(y,r),z) \; dz-\int_{\er^2} \deg(g,B(y,r),z) \; dz \Bigr|<\epsilon.\\
\end{aligned}
$$
Indeed, if this were not true, we would have uniformly
converging sequence such that conclusion of \eqref{reallyclose} would not hold. Analogously to the proof of Lemma \ref{lemma:DegreeLimit} we would then get a contradiction.

Moreover, similarly to the proof of Lemma \ref{lemma:DegreeLimit}, the Sobolev
inequality gives for these a.e.\ radii
  \begin{align}\label{hei}
    \Bigl|\int_{\er^2} \deg(f_s,B(y,r),z)\; dz\Bigr|
    \leq C | D f_{y,s} | (\partial B(y,s)).
  \end{align}

Given $\epsilon>0$ we choose $\eta>0$ as in \eqref{reallyclose} and then we choose $\delta>0$ so that for every
$s\in[r-\delta,r]$ we have
\begin{equation}
\label{laskukiitos}
\|f-f_s\|_{L^{\infty}(\partial B)}<\eta\ \text{ and }\ \frac{1}{\delta}\int_{r -\delta}^r\Bigl||Df_{y,s}|(\partial
B(y,s))-|Df_{y,r}|(\partial B(y,r))\Bigr|\; ds<\eta^2
\end{equation}
where we have used \eqref{lebpoint}. By Chebyshev's inequality with
\eqref{laskukiitos} we obtain
$$
|W|<\eta\delta\text{ for }W:=\Bigl\{s\in [r-\delta,r]:\
\bigl||Df_{y,s}|(\partial
B(y,s))-|Df_{y,r}|(\partial B(y,r))\bigr|>\eta\Bigr\}.
$$

By \eqref{ahoj}, \eqref{444}, $\int_{r-\delta}^r \Phi'_{\delta}=1$,
$|\Phi'_{\delta}|\leq \frac{C}{\delta}$, \eqref{reallyclose} and \eqref{hei} %and
%\eqref{tauchoice}
we obtain
$$
\begin{aligned}
\Bigl|\int_{\Omega}  &\Phi_\delta(|x-y|)\; d\mathcal J_f(x)- \int_{\er^2} \deg(f,B(y,r),z) \; dz\Bigr|=\\
&=\Bigl| \int_{r-\delta}^{r} \Phi'_{\delta}(s) \Bigl(\int_{\er^2} \deg (f_s, B(y,r), z)\;  dz- \int_{\er^2} \deg (f, B(y,r), z)\;  dz\Bigr)\; ds\Bigr|\\
&\leq \frac{C}{\delta} \Bigl[\int_{[r-\delta,r]\setminus W}\epsilon
+\int_W \bigl(| D f_{y,s} |(\partial B(y,s))+| D f_{y,r}| (\partial B(y,r))\bigr)\; ds\Bigr]\\
&\leq C\epsilon+\frac{C}{\delta} \int_W \bigl|| D f_{y,s} |(\partial B(y,s))-| D
f_{y,r}| (\partial B(y,r))\bigr|\; ds+ \frac{2 C}{\delta}\int_W |Df_{y,r}|(
\partial B(y,r))ds\\
&\leq C\epsilon+C\eta^2+2 C\eta |Df_{y,r}|( \partial B(y,r)).\\
\end{aligned}
$$

  Together with \eqref{eq:InnerRegularity} this implies that
  $$
    \mathcal J_f(B(y,r))
    = \lim_{\delta\to 0+} \int_\Omega \Phi_\delta(|x-y|)\; d\mathcal J_f(x)
    =\int_{\er^2} \deg (f, B(y,r), z)\;  dz.
  $$
\end{proof}

\section{Proof of main Theorem \ref{regularity}}
Given a point $x\in\er^n$ and a $s>0$ we denote by $Q(x,s)$ the cube with center $x$ and sidelength $s$ and whose sides are parallel to coordinate planes. Given a $t>0$ we also denote $tQ(x,s)\colonequals Q(x,ts).$
\begin{proof}[Proof of Theorem \ref{regularity}]
  Without loss of  generality we may assume that $(-1,2)^3\subset\Omega$ and we
prove only that $f^{-1}\in BV\bigl(f((0,1)^3)\bigr)$  as the statement is local. We denote $Q\colonequals Q((\frac12,\frac12), 1)= (0,1)^2.$ Slightly abusing the notation we write $2Q\colonequals Q((\frac12,\frac12), 2)$.
  We claim that
  $$
  \int_{0}^1 \haus^2\bigl(f(Q\times \{t\})\bigr) dt<\infty
  $$
  and the statement of the theorem then follows from Theorem \ref{abstract}. Let $\epsilon>0$.
  We start with an estimate for $\haus^2_{\epsilon}(f(Q\times \{t\}))$ for some fixed $t\in(0,1).$

  First let us fix $t\in(0,1)$  such that (see \eqref{stupid})
	$$
	A(f,Q\times\{t\})<\infty
	$$
	and
\eqn{star}
  $$
  \lim_{\delta\to 0}\frac{1}{2\delta}\int_{t-\delta}^{t+\delta}
  \sum_{j=1}^3\bigl||\djac_{f^s_{3,j}}|(2Q\times\{s\})-|\djac_{f^t_{3,j}}|(2Q\times\{t\})\bigr|\; ds=0
  $$
  and we note that this holds for a.e.\ $t\in(0,1)$ by the Lebesgue density theorem.
  Let us define the measure on $(0,1)$ by
  \eqn{defmu}
  $$
  \begin{aligned}
    \mu((a,b))=&\sum_{j=1}^3\int_{-1}^2 |\djac_{f^s_{2,j}}|\bigl((-1,2)\times\{s\}\times(a,b)\bigr)\; ds \\
    &+\sum_{j=1}^3\int_{-1}^2 |\djac_{f^s_{1,j}}|\bigl(\{s\}\times(-1,2)\times(a,b)\bigr)\; ds. \\
  \end{aligned}
  $$
  Let us denote by $h$ the absolutely continuous part of $\mu$ with respect to $\mathcal{L}^1$. Then it is easy to see that
  \eqn{easy}
  $$
  \int_0^1 h\leq \mu((0,1))\leq |\dadj Df|\bigl((-1,2)^2\times (0,1)\bigr).
  $$
  Moreover, we can fix $t$ so that
  $$
  \lim_{\delta \to 0}\frac{\mu\bigl((t-\delta,t+\delta)\bigr)}{2\delta}=h(t)
  $$
  which holds for a.e.\ $t$ by the Lebesgue density theorem and by the fact
  that the corresponding limit is zero a.e.\ for the singular part of $\mu$.

  Since $f$ is uniformly continuous there exists
  for our fixed $t$ a subdivision of $Q\times\{t\}=\bigcup_i Q_i $ into squares $Q_i=Q(c^i,r_i)$
  such that $\diam(f(2Q_i\times \{t\}))< \frac{\epsilon}{2}$.
  Furthermore we fix $r>0$ so that for every $0<\delta<r$ we have with the help of \eqref{star}
  \eqn{ahojky}
  $$
  \begin{aligned}
    (i)\ &\frac{\mu\bigl((t-\delta,t+\delta)\bigr)}{2\delta}\leq 2 h(t),\\
    (ii)\ &\sum_{j=1}^3|\djac_{f^{t+\delta}_{3,j}}|(2Q\times\{t+\delta\})\leq 2 \sum_{j=1}^3|\djac_{f^{t}_{3,j}}|(2Q\times\{t\}),\\
    (iii)\ &\sum_{j=1}^3|\djac_{f^{t-\delta}_{3,j}}|(2Q\times\{t-\delta\})\leq 2 \sum_{j=1}^3|\djac_{f^{t}_{3,j}}|(2Q\times\{t\})\text{ and }\\
    (iv)\ &\diam\bigl(f(2Q_i\times [t-\delta,t+\delta])\bigr)< \epsilon\text{ for
each }i.
  \end{aligned}
  $$

  For $\eta>0$ we put our $Q_i\times \{t\}$ into the box
  $$
  U_{i,\eta}\colonequals (1+\eta)Q_i\times[t-\delta,t+\delta].
  $$
  In the following we divide $\partial  U_{i,\eta}$ into three parts parallel to
coordinate axes
  $$
  \partial_3 U_{i,\eta}\colonequals  (1+\eta)Q_i\times\{t-\delta,t+\delta\},\
\partial_2 U_{i,\eta}\text{ and }\partial_1 U_{i,\eta},
  $$
  where $\partial_2 U_{i,\eta}$  denotes two rectangles perpendicular to $x_2$
axis and $\partial_1 U_{i,\eta}$ denotes two rectangles perpendicular to $x_1$
axis.
  For each $Q_i$ we choose a real number $0\leq \eta_{i}\leq 1$ so that
  \begin{equation}
    \label{etachoice}
    \sum_{k=1}^2\sum_{j=1}^3|\djac_{f_{k,j}}|(\partial_k U_{i,\eta_{i}}) \leq
    \int_0^1 \sum_{k=1}^2\sum_{j=1}^3 |\djac_{f_{k,j}}| (\partial_k U_{i,\eta'}) \; d\eta',
  \end{equation}
  which is possible as the smallest value is less or  equal to the average and
here and in the following we denote for simplicity $|\djac_{f_{1,j}}|(\partial_1
U_{i,\eta})$ the sum of two
  $$
  |\djac_{f^{c^i_1\pm  (1+\eta)r_i}_{1,j}}|\bigl(\{c^i_1\pm (1+\eta)r_i\}\times[c^i_2-
(1+\eta)r_ i,c^i_2+ (1+\eta)r_i]\times[t-\delta,t+\delta] \bigr),
  $$
  as $Q_i=Q(c^i,r_ i)$.

  It is obvious that $f(Q_i\times\{t\})\subset f(U_{i,\eta_i})$.
  By the definition of  the Lebesgue area \eqref{lebarea} and its estimate
\eqref{Leb} we obtain that we can approximate $f$ on $\partial
U_{i,\eta_i}$ by piecewise linear $f^i\colon\partial U_{i,\eta_i}\to\er^3$ such
that
\eqn{area}
$$
    \haus^2(f^i(\partial U_{i,\eta_i}))
    \leq 2 A(f,\partial U_{i,\eta_i})
    =2\sum_{k=1}^3  A(f_k,\partial_k U_{i,\eta_i})\leq 2 \sum_{k=1}^3\sum_{j=1}^3
V(f_{k,j},\partial_k U_{i,\eta_i})
  $$
  and so that $f^i$ is so close to $f$ that (see \eqref{ahojky} $(iv)$)
  \eqn{ttt}
  $$
  \diam\bigl(f^i(\partial U_{i,\eta_i})\bigr)< \epsilon
  $$
  and $f(Q_i\times\{t\})$ lies inside $f^i(U_{i,\eta_i})$, i.e.
  $$
  f(Q_i\times\{t\})\subset G_i\colonequals \bigcup \text{ bounded components of }\er^3\setminus f^i(\partial U_{i,\eta_i}).
  $$
  By \eqref{ttt} we have
  $$
  \haus^2_\epsilon(G_i)=\haus^2_\infty(G_i).
  $$

  Now we obtain for each $t$ with the help of Gustin boxing lemma \eqref{gustin} and \eqref{area}
  \begin{align}
    \label{splitting}
    \nonumber\haus^2_\epsilon(f(Q\times\{t\}))
    &\leq \sum_i \haus^2_\epsilon(f(Q_i\times\{t\}))\leq \sum_i \haus^2_\epsilon(G_i)\\
    &=\sum_i \haus^2_\infty(G_i)
      \leq C \sum_i \haus^2(\partial G_i)\leq C \sum_i \haus^2( f^i(\partial U_{i,\eta_{i}}))\\
    &\leq\nonumber
      C \sum_i\sum_{k=1}^3\sum_{j=1}^3 V(f_{k,j},\partial_k U_{i,\eta_i}) .
  \end{align}
  Recall that by definition \eqref{recalldef}
  \begin{align*}
    V(h,U)
    = \sup_S\left\{\sum_{\pi\in S}\int_{\er^2} \bigl|\deg(h,\pi,y)\bigr|\;
dy\right\},
  \end{align*}
  so by Proposition \ref{coro:degreeWithAbsoluteValues} we obtain
        \eqn{star2}
  $$
  \begin{aligned}
    \sum_{k=1}^2\sum_{j=1}^3 V(f_{k,j},  \partial_{k} U_{i,\eta_{i}})&\leq C
\sum_{k=1}^2\sum_{j=1}^3|\djac_{f_{k,j}}|(\partial_{k}U_{i,\eta_{i}}) \text{ and
}\\
    \sum_{j=1}^3 V(f_{3,j}, \partial_{3}  U_{i,\eta_{i}})&\leq C
\sum_{j=1}^3|\djac_{f_{3,j}}|(\partial_{3}U_{i,\eta_{i}}).\\
  \end{aligned}
  $$
Notice that even though Theorem \ref{coro:degreeWithAbsoluteValues} is stated
only for disks, it also holds for rectangles and moreover, we may use it for polygons.
This can be seen by covering the
polygon by rectangles and arguing as in the end of the
proof of Proposition \ref{coro:degreeWithAbsoluteValues}.

  We treat the terms in \eqref{star2} separately. We sum  the last inequality, use the fact that
$(1+\eta_i)Q_i$ have bounded overlap (as $1\leq 1+\eta_i\leq 2$) and with the
help of \eqref{ahojky} $(ii)$ and $(iii)$ we obtain
  \eqn{hor}
  $$
  \begin{aligned}
    \sum_i\sum_{j=1}^3 V(f_{3,j}, &\partial_{3} U_{i,\eta_{i}})
    \leq C \sum_i\sum_{j=1}^3|\djac_{f_{3,j}}|(\partial_{3}U_{i,\eta_{i}})\\
    &\leq C
\sum_{j=1}^3\Bigl(|\djac_{f_{3,j}}|(2Q\times\{t-\delta\})+|\djac_{f_{3,j}}
|(2Q\times\{t+\delta\}\Bigr)\\
    &\leq C \sum_{j=1}^3 |\djac_{f_{3,j}}|(2Q\times\{t\}).
  \end{aligned}
  $$
  For the remaining part of the right  hand side of \eqref{splitting} we recall that
$Q_i=Q(c^i,r_i)$ and by \eqref{star2}, \eqref{etachoice}, linear change of variables and $\delta<r$
we have
  $$
  \begin{aligned}
    \sum_{k=1}^2\sum_{j=1}^3 V(f_{k,j},& \partial_{k} U_{i,\eta_{i}})
    \leq C \sum_{k=1}^2\sum_{j=1}^3|\djac_{f_{k,j}}|(\partial_{k}U_{i,\eta_{i}}) \\
    &\leq C\int_0^1 \sum_{k=1}^2\sum_{j=1}^3 |\djac_{f_{k,j}}| (\partial_k U_{i,\eta'}) \; d\eta'\\
    &\leq \frac{C}{\delta} \int_{c^i_1-2r_i}^{c^i_1+2r_i}
    \sum_{j=1}^3 |\djac_{f^a_{1,j}}|\bigl(\{a\}\times[c^i_2-2 r_i,c^i_2+2 r_i]\times[t-\delta,t+\delta]\bigr)\; da\\
    &+\frac{C}{\delta} \int_{c^i_2-2r_i}^{c^i_2+2r_i}
    \sum_{j=1}^3 |\djac_{f^a_{2,j}}|\bigl([c^i_1-2 r_i,c^i_1+2 r_i]\times\{a\}\times[t-\delta,t+\delta]\bigr)\; da.
  \end{aligned}
  $$
  Summing over $i,$ using bounded overlap of $2 Q_i$, \eqref{defmu} and \eqref{ahojky} $(i)$  we obtain
  \eqn{vert}
  $$
  \begin{aligned}
    \sum_i \sum_{k=1}^2\sum_{j=1}^3 V(f_{k,j}, \partial_{k} U_{i,\eta_{i}})
    \leq C \frac{\mu\bigl((t-\delta,t+\delta) \bigr)}{\delta}\leq C h(t).\\
  \end{aligned}
  $$

  Combining \eqref{splitting}, \eqref{hor} and \eqref{vert}, we have with the help of \eqref{easy}
  $$
  \begin{aligned}
    \int_0^1\haus^2_\epsilon\bigl(f(Q\times \{t\})\bigr)\; dt
    &\leq C\int_0^1\sum_{j=1}^3 |\djac_{f_{3,j}}|(2Q\times\{t\})\; dt+\int_0^1 h(t)\; dt\\
    &\leq C   |\dadj Df|((-1,2)^3).\\
  \end{aligned}
  $$
  By passing $\epsilon\to 0$ we obtain our conclusion with the help of Theorem \ref{abstract}.
\end{proof}

\section{Reverse implication}

The main aim of this Section is to show Theorem \ref{reverse}. For its proof we again use some ideas from M\"uller \cite{M} and De Lellis \cite{D}.
As a corollary we show that the notion of
$\dadj Df\in \mathcal M$ does not depend on the chosen system of coordinates and
that this notion is weakly closed.

For the proof of Theorem \ref{reverse} we require the following result
which shows that the topological degree is smaller than the number of preimages.
\prt{Lemma}
\begin{proclaim}\label{lemma}
 Let $F \colon \er^3 \to \er^3$ be a homeomorphism,
 $f \colon \er^2 \to \er^3$ the restriction of $F$ to the $xy$-hyperplane,
 $p \colon \er^3 \to \er^2$ the projection $(x_1,x_2,x_3) \mapsto (x_1,x_2)$
 and $g = p \circ f$.
 Then for any $B(x,r) \subset (0,1)^2$ and $y \in \er^2 \setminus g (\partial B(z,r))$,
 $$
 |\deg (g, B(z,r), y)|\leq  N(g, B(z,r), y).
 $$
\end{proclaim}
\begin{proof}
  We may assume that $N(g, B(z,r), y)$ is finite and
  $\deg(g, B(z,r),y) > 0$.
  By \cite[Theorem 2.9]{FG} we can write the degree as a sum of local indices
  $$
  \deg(g, B(z,r), y)
  =\sum_{x \in B(z,r) \cap g\inv\{y\}}i(g,x,y);
  $$
  recall that local index is defined by
  $$
  i(g,x,y):=\deg(g,V,y),
  $$
  where $V$ is any neighborhood of $x$ such that $g^{-1}\{y\}\cup\bar{V} =\{ x \}$.

  To prove the claim it thus suffices to prove that $|i(g,x,y)|\leq 1$ for every
  $x\in g \inv \{ y \} $. Towards contradiction suppose that this is not the case.
  Fix some $x_0\in g \inv \{ y \}$ such that $i(g,x,y)\geq 2$; the case when
  the index is negative is dealt identically. Let $B(x_0,s)$ be a ball such that
  $$
  i(g,x_0,y)=\deg(g,B(x_0,s),y).
  $$

  Without loss of generality we may assume that $x_0 = y = 0$, $s=1$.
  Denote $Z = \{ 0 \} \times \{ 0 \} \times \er$.
  Since the topological degree equals the winding number, the
  path $\beta \colonequals g (\partial B(0,1))$ winds around the point $0$ at least twice
  in $\er^2 \setminus \{ 0 \}$, so especially the path
  $\alpha \colonequals f (\partial B(0,1) )$ winds twice around $Z$ in
  $\er^3 \setminus Z$.

  Now we note that $\partial B(0,1) \times \{ 0 \} \subset \es^2$, where $\es^2$ denotes the two-dimensional sphere in $\er^3$.
  Since
  $F\colon \er^3 \to \er^3$ is a a homeomorphism and $f$ the restriction of $F$,
  $F(\es^2)$ is a topological sphere
  in $\er^3$. %; without loss of generality we may assume that $F(\es^2) = \es^2$.
  Furthermore, $Z$ intersects $f B(0,1)$ only at a single point, and
  we fix $\hat Z$ to be the compact subinterval of $Z$ which contains
  the intersection point and intersects $F(\es^2)$ only at the endpoints of
  the interval, which we may assume to be $(0,0,\pm 1)$. The unique pre-images
  of these points cannot be on the circle $\partial B(0,1)$, so we may assume then
  to be $(0,0,\pm 1)$ as well. Thus
  $$
  \alpha
  = f (\partial B(0,1) )
  = F (\partial B(0,1) \times \{ 0 \})
  \subset F(\es^2) \setminus \hat Z.
  $$
  This gives rise to a contradiction, since $F$ is a homeomorphism and
  so the degree of $F|_{\es^2 \setminus Z}$
  is $\pm 1$. More specifically, the path $\alpha \colon \es^1 \to F(\es^2) \setminus Z$
  winds around the $Z$-axis at least twice, i.e.\
  the homotopy class $[\alpha]$ of $\alpha$ in the group
  $\pi_1(\er^2 \setminus Z, \alpha(0)) \simeq \zet$ is non-zero and does not
  span the group $\zet$.
  Furthermore the intersection $Z \cap F(B^3(0,1))$ consists of countably many paths starting
  and ending at the boundary $f\es^2$ and so
  since $Z$ intersects $f B(0,1)$ only at a single point all but one of these loops
  can be pulled to the boundary $f \es^2$ without intersecting $\alpha$.
  Thus the homotopy class $[\alpha]$ of $\alpha$ in the group
  $\pi_1(F(\es^2) \setminus \hat Z, \alpha(0)) \simeq \zet$ is also non-zero and does not
  span the group $\zet$. But this is a contradiction since
  $\alpha = F (\partial B(0,1) \times \{ 0 \})$, where the homotopy class $[ \partial B(0,1) \times \{ 0 \}]$
  spans $\pi_1(F(\es^2) \setminus Z, (1,0,0)) \simeq \zet$
  at the domain side and a homeomorphism $F$ induces an isomorphism between
  homotopy groups by e.g.\ \cite[p. 34]{Hatcher}.
\end{proof}

\begin{proof}[Proof of Theorem \ref{reverse}]
  The distributional adjugate is a well-defined  distribution as $f\in BV$ is
continuous. Without loss of generality we assume that $f$ is defined on
$(0,1)^3$ and we show that $\dadj Df\in\M((0,1)^3)$.
  We only show that $\djac_{f^t_{1,1}}$ is a measure for a.e.\ $t\in(0,1)$ and that
  \eqn{b}
  $$
  \int_0^1 \djac_{f^t_{1,1}}((0,1)^2)\; dt<\infty
  $$
  as the proof for other eight components of $\dadj Df$ is similar.

  By Theorem \ref{abstract} we know that
  $$
  \int_0^1 \haus^2\bigl( f^t_{1}((0,1)^2) \bigr)\;
dt<\infty
  $$
   and hence
	$$
	A(f,[0,1]^2\times\{t\})<\infty \text{ for a.e. }t\in(0,1). 
	$$
  Let us fix $t\in (0,1)$ such that $\haus^2\bigl( f^t_{1}((0,1)^2) \bigr)<\infty$.
  Put $g \colonequals f^t_{1,1}$ and denote by $g_1$ and $g_2$ its coordinate functions.
  Let us fix $\varphi \in C^{1}_0((0,1)^2)$. We recall the definition of distributional Jacobian
  \begin{align*}
    \mathcal J_g(\varphi)
    &=-\int_{(0,1)^2}g_1(x) J(\varphi(x),g_2(x))\; dx \\
    &=-\int_{(0,1)^2} \Bigl\langle [g_1(x),0] \cdot \cof Dg(x),D\varphi(x) \Bigr\rangle\; dx,
  \end{align*}
  where the integration is with respect to the relevant components of the
variation measure of $g$ as earlier.

  Let $\psi\in C^{\infty}_C[0,1)$ be such that $\psi\geq 0$, $\psi' \leq 0$ and
  \begin{align*}
    \int_{B(0,1)}\psi(|x|)\;dx
    =1.
  \end{align*}
  For each $\eps>0$ we denote by $\eta_\eps$ the usual convolution kernel, that
is
  \begin{align*}
    \eta_\eps(x)
    = {\psi_\eps(|x|)}
    = \eps^{-2}\psi \Bigl(\frac{|x|}{\eps}\Bigr).
  \end{align*}
  It is clear that $\eta_\eps*D\phi=D\eta_\eps*\phi$ converges uniformly to
$D\phi$ as $\eps\to 0+$ and hence
  $$
  \mathcal J_g(\varphi)=\lim_{\eps\rightarrow 0_+} - \int_{(0,1)^2} \Bigl\langle
[g_1(x),0] \cdot \cof Dg(x),
  \Bigl ( \int_{B(x,\eps)}\varphi D \eta_\eps (x-z)\; dz\Bigr)\Bigr \rangle\;
dx.
  $$
  It is easy to see that $D\eta_\eps(x)={\psi_\eps}'(|x|)\nu$, where
$\nu=\frac{x}{|x|}$ is the normal vector.
  By the Fubini theorem and change to polar coordinates we get
$$
    \mathcal J_g(\varphi)=
    -\lim_{\eps\rightarrow 0_+}\int_{(0,1)^2}\varphi (z)\Bigl ( \int_0^\eps
\psi'_\eps(r)\int_{\partial B(z,r)} g_1(x) d(Dg|_{\partial B(z,r)}(x)) dr\Bigr)
dz.
$$
  By the degree formula Proposition \ref{prop:degree} we obtain
  \begin{align*}
    \mathcal J_g(\varphi)= -\lim_{\eps\rightarrow 0_+}\int_{(0,1)^2} \varphi(z)
    \Bigl ( \int_0^\eps \psi'_\eps(r)\int_{\er^2} \deg (g, B(z,r), y)\; dy\; dr
\Bigr) dz.
  \end{align*}
  Let us fix $0<\eps<\frac12\dist(\supp\,\varphi,\partial(0,1)^2)$. Then we have with the help of Lemma \ref{lemma}
  \eqn{tere}
  $$
  \begin{aligned}
    |\mathcal J_g(\varphi)|&
    \leq 2 \int_{\spt(\phi)} |\varphi(z)|\Bigl ( \int_0^\eps
|\psi'_\eps(r)|\int_{\er^2} |\deg (g, B(z,r), y)|\; dy\; dr \Bigr) dz\\
    &\leq 2\|\phi\|_{\infty}\int_{(0,1)^2} \Bigl (
\int_0^\eps \frac{C}{\eps^3}\int_{\er^2} N(g, B(z,r), y)\; dy\; dr \Bigr) dz\\
    &\leq C \|\phi\|_{\infty}\int_{\er^2}\frac{1}{\eps^2}
\int_{(0,1)^2} N(g, B(z,\eps), y)\; dz\; dy.
  \end{aligned}
  $$
Notice that for fixed $y\in\er^2$  we have
$$
N(g,B(z,t),y)=\sum_{z_i\in g^{-1}\{y\}}\chi_{B(z_i,t)}(z).
$$
With this we obtain from \eqref{tere} that
\eqn{teretaas}
$$
|\mathcal J_g(\varphi)|\leq C \|\phi\|_{\infty}\int_{\er^2} N(g,
(0,1)^2,y)\;
dy.
$$

  By \cite[Theorem 7.7]{Mattila} we see that
  \begin{align*}
    \int_{\er^2} N(g, (0,1)^2, y)\; dy
    \leq \haus^2( f^t_1((0,1)^2)).
  \end{align*}
  Combining this with \eqref{teretaas},
  it follows that for every $\phi\in C^{1}_0((0,1)^2)$ we have
  \eqn{a}
  $$
  |\djac_g(\phi)| \leq C \|\phi\|_{\infty}
\haus^2\bigl(f^t_{1}((0,1)^2)\bigr)
  $$
  with $C$ independent of $\varphi.$
  By the Hahn-Banach Theorem there is an extension to every $\phi\in C_0((0,1)^2)$ which satisfies the same bound. By the Riesz Representation Theorem there is a measure $\mu_t$ such that
  $$
  \djac_g(\phi)=\int_{(0,1)^2}\phi(x)\; d\mu_t(x)\text{ for every }\phi\in C^{1}_0((0,1)^2).
  $$
  By \eqref{a} and \eqref{b} we have
  $$
  \int_0^1 \mu_t((0,1)^2)\; dt\leq C\int_0^1 \haus^2\bigl( f^t_{1}((0,1)^2) \bigr)\; dt<\infty
  $$
  and thus $\dadj Df\in \M((0,1)^3)$.
\end{proof}

\subsection{Dependence on the system of coordinates}
\label{sec:ADJ-CoordinateDependence}

In principle the Definition \ref{def} of $\dadj Df\in\M$ depends on our coordinate system. Below we show that this notion is independent on the system of coordinates.

\prt{Corollary}
\begin{proclaim}
  Let $\Omega\subset\er^3$ be a domain and $f\in BV(\Omega,\er^3)$ be a
  homeomorphism satisfying \eqref{stupid} such that $\dadj Df\in \M(\Omega,\er^{3\times 3})$. Then $\dadj Df\in \M(\Omega,\er^{3\times 3})$ also for a different coordinate system.
\end{proclaim}
\begin{proof}
  By Theorem \ref{regularity} we know that $f^{-1}\in BV$. Hence $f\in BV_{\loc}$ and $f^{-1}\in BV_{\loc}$ and both of these do not depend on the choice of coordinate system. Thus by Theorem \ref{reverse} we have $\dadj Df\in \M(\Omega,\er^{3\times 3})$ for any coordinate system.
\end{proof}

It is of course not true that the value of
$$
|\dadj Df|(\Omega)
$$
is independent of coordinate system. In fact it might be more natural to define $|\dadj Df|$ as an average over all directions (and not only 3 coordinate directions). Then, one could ask for the validity of (compare
with \eqref{identity})
$$
|Df^{-1}|(f(\Omega))=|\dadj Df|(\Omega).
$$

\subsection{The notion is stable under weak convergence}

For possible applications in the Calculus of Variations we need to know that the notion of distributional adjugate is stable under weak convergence.

\prt{Theorem}
\begin{proclaim}
  Let $\Omega\subset\er^3$ be a bounded domain.
  Let $f_j, f$ be a $BV$ homeomorphisms of $(0,1)^3$ onto $\Omega$ and assume that $f_j\to f$ uniformly and weak* in $BV((0,1)^3,\Omega)$. Further suppose each $f_j$ satisfies \eqref{stupid} and  let $\dadj Df_j\in\M((0,1)^3)$ with
  \eqn{ppp}
  $$
  \sup_j |\dadj Df_j|\bigl((0,1)^3\bigr)<\infty.
  $$
  Then $\dadj Df\in\M$.
\end{proclaim}
\begin{proof}
  By \eqref{ppp} and Theorem \ref{regularity} we obtain that the sequence $(f^{-1}_j)$ is a bounded in $BV(\Omega,\er^3)$ and hence it has a weakly* converging subsequence. Thus we can assume (passing to a subsequence) that  $f^{-1}_j\to h$ weakly* in $BV$ and also strongly in $L^1$ (see \cite[Corollary 3.49]{AFP}). We define the pointwise representative of $h$ as
  $$
  h(y)\colonequals \limsup_{r\to 0}\frac{1}{|B(y,r)|}\int_{B(y,r)}h.
  $$

  Now we need to show that $h=f^{-1}$. Fix $x_0\in(0,1)^3$ and $0<r<\dist(x_0,\partial (0,1)^3)$. We find $\delta>0$ so that $B(f(x_0),\delta)$ is compactly contained in $f(B(x_0,r))$. Since $f_j \to f$ uniformly we obtain that for $j$ large enough we have
  $$
  B(f(x_0),\delta)\subset f_j(B(x_0,r)).
  $$
  It follows that
  $$
  f^{-1}_j(B(f(x_0),\delta))\subset B(x_0,r)\text{ and hence }|h(f(x_0))-x_0|\leq r
  $$
  where we use that $f_j^{-1}\to h$ strongly in $L^1$ and that we have a proper representative of $h$. As the above inequality holds for every $r>0$ we obtain $h(f(x_0))=x_0$.

  From $f\in BV$ and $f^{-1}=h\in BV$ we obtain $\dadj Df\in\M((0,1)^3)$ by Theorem \ref{reverse}.
\end{proof}

\vskip 5pt
\noindent
{\bf Acknowledgments.} The authors would like to thank Jan Mal\'y for pointing their interest to Theorem \ref{coarea2} and for many valuable comments  and for finding the gap in the original proof of Proposition \ref{coro:degreeWithAbsoluteValues}. The authors would also like to thank to Ulrich Menne for his information about literature on Lebesgue area and the anonymous referee for careful reading of the manuscript.

\end{document}